%
%
\def\Date{2008/05/09}


\ifx\pdfoutput\jamaisdefined\else
\input supp-pdf.tex \pdfoutput=1 \pdfcompresslevel=9

\fi

%

\magnification=1200
\hsize=11.25cm
\vsize=18cm
\parskip 0pt
\parindent=12pt
\voffset=1cm
\hoffset=1cm



\catcode'32=9

\font\tenpc=cmcsc10
\font\eightpc=cmcsc8
\font\eightrm=cmr8
\font\eighti=cmmi8
\font\eightsy=cmsy8
\font\eightbf=cmbx8
\font\eighttt=cmtt8
\font\eightit=cmti8
\font\eightsl=cmsl8
\font\sixrm=cmr6
\font\sixi=cmmi6
\font\sixsy=cmsy6
\font\sixbf=cmbx6

\skewchar\eighti='177 \skewchar\sixi='177
\skewchar\eightsy='60 \skewchar\sixsy='60

\catcode`@=11

\def\tenpoint{%
  \textfont0=\tenrm \scriptfont0=\sevenrm \scriptscriptfont0=\fiverm
  \def\rm{\fam\z@\tenrm}%
  \textfont1=\teni \scriptfont1=\seveni \scriptscriptfont1=\fivei
  \def\oldstyle{\fam\@ne\teni}%
  \textfont2=\tensy \scriptfont2=\sevensy \scriptscriptfont2=\fivesy
  \textfont\itfam=\tenit
  \def\it{\fam\itfam\tenit}%
  \textfont\slfam=\tensl
  \def\sl{\fam\slfam\tensl}%
  \textfont\bffam=\tenbf \scriptfont\bffam=\sevenbf
  \scriptscriptfont\bffam=\fivebf
  \def\bf{\fam\bffam\tenbf}%
  \textfont\ttfam=\tentt
  \def\tt{\fam\ttfam\tentt}%
  \abovedisplayskip=12pt plus 3pt minus 9pt
  \abovedisplayshortskip=0pt plus 3pt
  \belowdisplayskip=12pt plus 3pt minus 9pt
  \belowdisplayshortskip=7pt plus 3pt minus 4pt
  \smallskipamount=3pt plus 1pt minus 1pt
  \medskipamount=6pt plus 2pt minus 2pt
  \bigskipamount=12pt plus 4pt minus 4pt
  \normalbaselineskip=12pt
  \setbox\strutbox=\hbox{\vrule height8.5pt depth3.5pt width0pt}%
  \let\bigf@ntpc=\tenrm \let\smallf@ntpc=\sevenrm
  \let\petcap=\tenpc
  \normalbaselines\rm}

\def\eightpoint{%
  \textfont0=\eightrm \scriptfont0=\sixrm \scriptscriptfont0=\fiverm
  \def\rm{\fam\z@\eightrm}%
  \textfont1=\eighti \scriptfont1=\sixi \scriptscriptfont1=\fivei
  \def\oldstyle{\fam\@ne\eighti}%
  \textfont2=\eightsy \scriptfont2=\sixsy \scriptscriptfont2=\fivesy
  \textfont\itfam=\eightit
  \def\it{\fam\itfam\eightit}%
  \textfont\slfam=\eightsl
  \def\sl{\fam\slfam\eightsl}%
  \textfont\bffam=\eightbf \scriptfont\bffam=\sixbf
  \scriptscriptfont\bffam=\fivebf
  \def\bf{\fam\bffam\eightbf}%
  \textfont\ttfam=\eighttt
  \def\tt{\fam\ttfam\eighttt}%
  \abovedisplayskip=9pt plus 2pt minus 6pt
  \abovedisplayshortskip=0pt plus 2pt
  \belowdisplayskip=9pt plus 2pt minus 6pt
  \belowdisplayshortskip=5pt plus 2pt minus 3pt
  \smallskipamount=2pt plus 1pt minus 1pt
  \medskipamount=4pt plus 2pt minus 1pt
  \bigskipamount=9pt plus 3pt minus 3pt
  \normalbaselineskip=9pt
  \setbox\strutbox=\hbox{\vrule height7pt depth2pt width0pt}%
  \let\bigf@ntpc=\eightrm \let\smallf@ntpc=\sixrm
  \let\petcap=\eightpc
  \normalbaselines\rm}
\catcode`@=12

\tenpoint


\long\def\irmaaddress{{%
\bigskip
\eightpoint
\rightline{\quad
\vtop{\halign{\hfil##\hfil\cr
I.R.M.A. UMR 7501\cr
Universit\'e Louis Pasteur et CNRS,\cr
7, rue Ren\'e-Descartes\cr
F-67084 Strasbourg, France\cr
{\tt guoniu@math.u-strasbg.fr}\cr}}\quad}
}}


\font\tengoth=eufm10
\def\goth#1{\hbox{\tengoth #1}}


\catcode`\@=11
\def\pc#1#2|{{\bigf@ntpc #1\penalty \@MM\hskip\z@skip\smallf@ntpc%
	\uppercase{#2}}}
\catcode`\@=12

\def\pointir{\discretionary{.}{}{.\kern.35em---\kern.7em}\nobreak
   \hskip 0em plus .3em minus .4em }

\def\qed{\quad\raise -2pt\hbox{\vrule\vbox to 10pt{\hrule width 4pt
   \vfill\hrule}\vrule}}

\def\rem#1|{\par\medskip{{\it #1}\pointir}}

\def\vspace[#1]{\noalign{\vskip#1}}

\def\abstract#1{\vbox{\eightpoint\narrower\narrower 
\pc ABSTRACT|\pointir #1}}


\def\section#1{\goodbreak\par\vskip .3cm\centerline{\bf #1}
   \par\nobreak\vskip 3pt }

\long\def\th#1|#2\endth{\par\medbreak
   {\petcap #1\pointir}{\it #2}\par\medbreak}

\def\article#1|#2|#3|#4|#5|#6|#7|
    {{\leftskip=7mm\noindent
     \hangindent=2mm\hangafter=1
     \llap{{\tt [#1]}\hskip.35em}{\petcap#2}\pointir
     #3, {\sl #4}, {\bf #5} ({\oldstyle #6}),
     pp.\nobreak\ #7.\par}}
\def\livre#1|#2|#3|#4|
    {{\leftskip=7mm\noindent
    \hangindent=2mm\hangafter=1
    \llap{{\tt [#1]}\hskip.35em}{\petcap#2}\pointir
    {\sl #3}, #4.\par}}
\def\divers#1|#2|#3|
    {{\leftskip=7mm\noindent
    \hangindent=2mm\hangafter=1
     \llap{{\tt [#1]}\hskip.35em}{\petcap#2}\pointir
     #3.\par}}



\catcode`\@=11
\def\c@rr@#1{\vbox{%
  \hrule height \ep@isseur%
   \hbox{\vrule width\ep@isseur\vbox to \t@ille{%
           \vfil\hbox  to \t@ille{\hfil#1\hfil}\vfil}%
            \vrule width\ep@isseur}%
      \hrule height \ep@isseur}}
\def\ytableau#1#2#3#4{\vbox{%
  \gdef\ep@isseur{#2}
   \gdef\t@ille{#1}
    \def\\##1{\c@rr@{$#3 ##1$}}
  \lineskiplimit=-30cm \baselineskip=\t@ille%
    \advance \baselineskip by \ep@isseur%
     \halign{%
      \hfil$##$\hfil&&\kern -\ep@isseur%
       \hfil$##$\hfil \crcr#4\crcr}}}%
\catcode`\@=12

\def\Grille{\setbox13=\vbox to 5mm{\hrule width 110mm\vfill}
\setbox13=\vbox{\offinterlineskip
   \copy13\copy13\copy13\copy13\copy13\copy13\copy13\copy13
   \copy13\copy13\copy13\copy13\box13\hrule width 110mm}
\setbox14=\hbox to 5mm{\vrule height 65mm\hfill}
\setbox14=\hbox{\copy14\copy14\copy14\copy14\copy14\copy14
   \copy14\copy14\copy14\copy14\copy14\copy14\copy14\copy14
   \copy14\copy14\copy14\copy14\copy14\copy14\copy14\copy14\box14}
\ht14=0pt\dp14=0pt\wd14=0pt
\setbox13=\vbox to 0pt{\vss\box13\offinterlineskip\box14}
\wd13=0pt\box13}


\def\fleche(#1,#2)\dir(#3,#4)\long#5{%
\noalign{\nointerlineskip\leftput(#1,#2){\vector(#3,#4){#5}}\nointerlineskip}}


\def\hfl#1#2#3{\smash{\mathop{\hbox to#3{\rightarrowfill}}\limits
^{\scriptstyle#1}_{\scriptstyle#2}}}

\def\gfl#1#2#3{\smash{\mathop{\hbox to#3{\leftarrowfill}}\limits
^{\scriptstyle#1}_{\scriptstyle#2}}}


 \message{`lline' & `vector' macros from LaTeX}
 \catcode`@=11
\def\{{\relax\ifmmode\lbrace\else$\lbrace$\fi}
\def\}{\relax\ifmmode\rbrace\else$\rbrace$\fi}
\def\newcount{\alloc@0\count\countdef\insc@unt}
\def\newdimen{\alloc@1\dimen\dimendef\insc@unt}
\def\newwrite{\alloc@7\write\chardef\sixt@@n}

\newwrite\@unused
\newcount\@tempcnta
\newcount\@tempcntb
\newdimen\@tempdima
\newdimen\@tempdimb
\newbox\@tempboxa

\def\@spaces{\space\space\space\space}
\def\@whilenoop#1{}
\def\@whiledim#1\do #2{\ifdim #1\relax#2\@iwhiledim{#1\relax#2}\fi}
\def\@iwhiledim#1{\ifdim #1\let\@nextwhile=\@iwhiledim
        \else\let\@nextwhile=\@whilenoop\fi\@nextwhile{#1}}
\def\@badlinearg{\@latexerr{Bad \string\line\space or \string\vector
   \space argument}}
\def\@latexerr#1#2{\begingroup
\edef\@tempc{#2}\expandafter\errhelp\expandafter{\@tempc}%
\def\@eha{Your command was ignored.
^^JType \space I <command> <return> \space to replace it
  with another command,^^Jor \space <return> \space to continue without it.}
\def\@ehb{You've lost some text. \space \@ehc}
\def\@ehc{Try typing \space <return>
  \space to proceed.^^JIf that doesn't work, type \space X <return> \space to
  quit.}
\def\@ehd{You're in trouble here.  \space\@ehc}

\typeout{LaTeX error. \space See LaTeX manual for explanation.^^J
 \space\@spaces\@spaces\@spaces Type \space H <return> \space for
 immediate help.}\errmessage{#1}\endgroup}
\def\typeout#1{{\let\protect\string\immediate\write\@unused{#1}}}

\font\tenln    = line10
\font\tenlnw   = linew10

\newdimen\@wholewidth
\newdimen\@halfwidth
\newdimen\unitlength 

\unitlength =1pt


\def\thinlines{\let\@linefnt\tenln \let\@circlefnt\tencirc
  \@wholewidth\fontdimen8\tenln \@halfwidth .5\@wholewidth}
\def\thicklines{\let\@linefnt\tenlnw \let\@circlefnt\tencircw
  \@wholewidth\fontdimen8\tenlnw \@halfwidth .5\@wholewidth}

\def\linethickness#1{\@wholewidth #1\relax \@halfwidth .5\@wholewidth}

\newif\if@negarg

\def\lline(#1,#2)#3{\@xarg #1\relax \@yarg #2\relax
\@linelen=#3\unitlength
\ifnum\@xarg =0 \@vline
  \else \ifnum\@yarg =0 \@hline \else \@sline\fi
\fi}

\def\@sline{\ifnum\@xarg< 0 \@negargtrue \@xarg -\@xarg \@yyarg -\@yarg
  \else \@negargfalse \@yyarg \@yarg \fi
\ifnum \@yyarg >0 \@tempcnta\@yyarg \else \@tempcnta -\@yyarg \fi
\ifnum\@tempcnta>6 \@badlinearg\@tempcnta0 \fi
\setbox\@linechar\hbox{\@linefnt\@getlinechar(\@xarg,\@yyarg)}%
\ifnum \@yarg >0 \let\@upordown\raise \@clnht\z@
   \else\let\@upordown\lower \@clnht \ht\@linechar\fi
\@clnwd=\wd\@linechar
\if@negarg \hskip -\wd\@linechar \def\@tempa{\hskip -2\wd\@linechar}\else
     \let\@tempa\relax \fi
\@whiledim \@clnwd <\@linelen \do
  {\@upordown\@clnht\copy\@linechar
   \@tempa
   \advance\@clnht \ht\@linechar
   \advance\@clnwd \wd\@linechar}%
\advance\@clnht -\ht\@linechar
\advance\@clnwd -\wd\@linechar
\@tempdima\@linelen\advance\@tempdima -\@clnwd
\@tempdimb\@tempdima\advance\@tempdimb -\wd\@linechar
\if@negarg \hskip -\@tempdimb \else \hskip \@tempdimb \fi
\multiply\@tempdima \@m
\@tempcnta \@tempdima \@tempdima \wd\@linechar \divide\@tempcnta \@tempdima
\@tempdima \ht\@linechar \multiply\@tempdima \@tempcnta
\divide\@tempdima \@m
\advance\@clnht \@tempdima
\ifdim \@linelen <\wd\@linechar
   \hskip \wd\@linechar
  \else\@upordown\@clnht\copy\@linechar\fi}

\def\@hline{\ifnum \@xarg <0 \hskip -\@linelen \fi
\vrule height \@halfwidth depth \@halfwidth width \@linelen
\ifnum \@xarg <0 \hskip -\@linelen \fi}

\def\@getlinechar(#1,#2){\@tempcnta#1\relax\multiply\@tempcnta 8
\advance\@tempcnta -9 \ifnum #2>0 \advance\@tempcnta #2\relax\else
\advance\@tempcnta -#2\relax\advance\@tempcnta 64 \fi
\char\@tempcnta}

\def\vector(#1,#2)#3{\@xarg #1\relax \@yarg #2\relax
\@linelen=#3\unitlength
\ifnum\@xarg =0 \@vvector
  \else \ifnum\@yarg =0 \@hvector \else \@svector\fi
\fi}

\def\@hvector{\@hline\hbox to 0pt{\@linefnt
\ifnum \@xarg <0 \@getlarrow(1,0)\hss\else
    \hss\@getrarrow(1,0)\fi}}

\def\@vvector{\ifnum \@yarg <0 \@downvector \else \@upvector \fi}

\def\@svector{\@sline
\@tempcnta\@yarg \ifnum\@tempcnta <0 \@tempcnta=-\@tempcnta\fi
\ifnum\@tempcnta <5
  \hskip -\wd\@linechar
  \@upordown\@clnht \hbox{\@linefnt  \if@negarg
  \@getlarrow(\@xarg,\@yyarg) \else \@getrarrow(\@xarg,\@yyarg) \fi}%
\else\@badlinearg\fi}

\def\@getlarrow(#1,#2){\ifnum #2 =\z@ \@tempcnta='33\else
\@tempcnta=#1\relax\multiply\@tempcnta \sixt@@n \advance\@tempcnta
-9 \@tempcntb=#2\relax\multiply\@tempcntb \tw@
\ifnum \@tempcntb >0 \advance\@tempcnta \@tempcntb\relax
\else\advance\@tempcnta -\@tempcntb\advance\@tempcnta 64
\fi\fi\char\@tempcnta}

\def\@getrarrow(#1,#2){\@tempcntb=#2\relax
\ifnum\@tempcntb < 0 \@tempcntb=-\@tempcntb\relax\fi
\ifcase \@tempcntb\relax \@tempcnta='55 \or
\ifnum #1<3 \@tempcnta=#1\relax\multiply\@tempcnta
24 \advance\@tempcnta -6 \else \ifnum #1=3 \@tempcnta=49
\else\@tempcnta=58 \fi\fi\or
\ifnum #1<3 \@tempcnta=#1\relax\multiply\@tempcnta
24 \advance\@tempcnta -3 \else \@tempcnta=51\fi\or
\@tempcnta=#1\relax\multiply\@tempcnta
\sixt@@n \advance\@tempcnta -\tw@ \else
\@tempcnta=#1\relax\multiply\@tempcnta
\sixt@@n \advance\@tempcnta 7 \fi\ifnum #2<0 \advance\@tempcnta 64 \fi
\char\@tempcnta}

\def\@vline{\ifnum \@yarg <0 \@downline \else \@upline\fi}

\def\@upline{\hbox to \z@{\hskip -\@halfwidth \vrule
  width \@wholewidth height \@linelen depth \z@\hss}}

\def\@downline{\hbox to \z@{\hskip -\@halfwidth \vrule
  width \@wholewidth height \z@ depth \@linelen \hss}}

\def\@upvector{\@upline\setbox\@tempboxa\hbox{\@linefnt\char'66}\raise
     \@linelen \hbox to\z@{\lower \ht\@tempboxa\box\@tempboxa\hss}}

\def\@downvector{\@downline\lower \@linelen
      \hbox to \z@{\@linefnt\char'77\hss}}

\thinlines

\newcount\@xarg
\newcount\@yarg
\newcount\@yyarg
\newcount\@multicnt
\newdimen\@xdim
\newdimen\@ydim
\newbox\@linechar
\newdimen\@linelen
\newdimen\@clnwd
\newdimen\@clnht
\newdimen\@dashdim
\newbox\@dashbox
\newcount\@dashcnt
 \catcode`@=12


\newbox\tbox
\newbox\tboxa

\def\leftzer#1{\setbox\tbox=\hbox to 0pt{#1\hss}%
     \ht\tbox=0pt \dp\tbox=0pt \box\tbox}

\def\rightzer#1{\setbox\tbox=\hbox to 0pt{\hss #1}%
     \ht\tbox=0pt \dp\tbox=0pt \box\tbox}

\def\centerzer#1{\setbox\tbox=\hbox to 0pt{\hss #1\hss}%
     \ht\tbox=0pt \dp\tbox=0pt \box\tbox}

%
\def\image(#1,#2)#3{\vbox to #1{\offinterlineskip
    \vss #3 \vskip #2}}


\def\leftput(#1,#2)#3{\setbox\tboxa=\hbox{%
    \kern #1\unitlength
    \raise #2\unitlength\hbox{\leftzer{#3}}}%
    \ht\tboxa=0pt \wd\tboxa=0pt \dp\tboxa=0pt\box\tboxa}

\def\rightput(#1,#2)#3{\setbox\tboxa=\hbox{%
    \kern #1\unitlength
    \raise #2\unitlength\hbox{\rightzer{#3}}}%
    \ht\tboxa=0pt \wd\tboxa=0pt \dp\tboxa=0pt\box\tboxa}

\def\centerput(#1,#2)#3{\setbox\tboxa=\hbox{%
    \kern #1\unitlength
    \raise #2\unitlength\hbox{\centerzer{#3}}}%
    \ht\tboxa=0pt \wd\tboxa=0pt \dp\tboxa=0pt\box\tboxa}

\unitlength=1mm

\def\cput(#1,#2)#3{\noalign{\nointerlineskip\centerput(#1,#2){#3}
                             \nointerlineskip}}


\ifx\pdfoutput\jamaisdefined
\input epsf

\fi


\parskip 0pt plus 1pt

\def\article#1|#2|#3|#4|#5|#6|#7|
    {{\leftskip=7mm\noindent
     \hangindent=2mm\hangafter=1
     \llap{{\tt [#1]}\hskip.35em}{#2},\quad %
     #3, {\sl #4}, {\bf #5} ({\oldstyle #6}),
     pp.\nobreak\ #7.\par}}
\def\livre#1|#2|#3|#4|
    {{\leftskip=7mm\noindent
    \hangindent=2mm\hangafter=1
    \llap{{\tt [#1]}\hskip.35em}{#2},\quad %
    {\sl #3}, #4.\par}}
\def\divers#1|#2|#3|
    {{\leftskip=7mm\noindent
    \hangindent=2mm\hangafter=1
     \llap{{\tt [#1]}\hskip.35em}{#2},\quad %
     #3.\par}}

\def\l{\lambda}

\def\mod{\mathop{\rm mod}}

\def\lig{\mathop{\goth g}}
\def\tmax{{\max_t}}

\def\enslettre#1{\font\zzzz=msbm10 \hbox{\zzzz #1}}
\def\setZ{\mathop{\enslettre Z}}

\def\setP{\mathop{\cal P}}

\def\setH{\mathop{\cal H}}
\def\setW{\mathop{\cal W}}


\rightline{\Date}
\bigskip

\centerline{\bf The Nekrasov-Okounkov hook length formula: }
\centerline{\bf refinement, elementary proof, extension and applications}
\bigskip
\centerline{Guo-Niu HAN}
\bigskip\medskip

\abstract{
The paper is devoted to the derivation of the expansion formula for the 
powers of the Euler Product in terms of partition hook lengths,
discovered by Nekrasov and Okounkov in their study of the Seiberg-Witten
Theory.  We provide a refinement based on a new property of $t$-cores,
and give an elementary proof by using the Macdonald identities.
We also obtain an extension by adding two more parameters, which appears 
to be a discrete interpolation between the Macdonald identities
and the generating function for $t$-cores. Several applications are 
derived, including the ``marked hook formula".
}

\bigskip

\def\itemm#1|#2| {{\leftskip=7mm\noindent \hangindent=2mm\hangafter=1 
    \llap{\S#1.\ }{#2}\par}}



\def\sec{1}
\section{\sec. Introduction} 
An explicit expansion formula for the powers of the Euler 
Product in terms 
of partition hook lengths was discovered by Nekrasov and Okounkov in 
their study of the Seiberg-Witten Theory [NO06] (see also [CO08], where
a Jack polynomial analogue was derived) and 
re-discovered by the author recently [Ha08a] using 
an appropriate hook length expansion technique [Ha08b]\footnote{$^*$}{%
\eightpoint
The author has indeed deposited a paper on arXiv [Ha08a; April, 2008]
that contained an explicit expansion formula for the powers of 
the Euler Product in terms of partition hook lengths.
A few days later he received an email from 
Andrei Okounkov who kindly pointed out that the expansion formula already 
appeared in his joint paper, which was deposited on arXiv 
in Section ``High Energy Physics - Theory" [NO06; June, 2003; 90 pages].
Although the ultimate formula is the same in both papers, the methods
of proof belong to different cultures. The author's original paper has
remained on arXiv.  The present one contains parts of it,
plus several new results.

}.  
In the present paper we again take up the study of the Nekrasov-Okounkov 
formula and 
obtain several results in the following four directions:

(1) We establish new properties on $t$-cores, which can be seen
as a refinement of the Nekrasov-Okounkov formula. 
The proof involves a bijection 
between $t$-cores and integer vectors constructed by Garvan, Kim and 
Stanton [GKS90].

(2) We provide an elementary proof of the  Nekrasov-Okounkov formula
by using the Macdonald identities for $A_\ell^{(a)}$ [Ma72]
and the properties on $t$-cores mentioned in (1).

(3) We obtain an extension 
by adding two more parameters~$t$ and $y$,
so that the resulting formula appears to be a  discrete interpolation 
between the Macdonald identities 
and the generating function for $t$-cores (see Corollary 5.3).
Our extension
opens the way to richer
specializations, including the generating function for partitions,
the Jacobi triple product identity,
the Macdonald identity for $A_\ell^{(a)}$, the classical hook length formula, 
the marked hook formula [Ha08a], the generating function for $t$-cores, and
the $t$-core analogues of the hook formula and of the marked hook
formula. 
We also prove another extension of the generating functions for $t$-cores.

(4)~As applications, we derive some new formulas about hook lengths, 
including the ``marked hook formula". We also improve a result 
due to Kostant [Ko04]. A hook length expression of integer value is obtained by
using the Lagrange inversion formula.

\medskip
The basic notions needed here can be found in 
[Ma95, p.1; St99, p.287;  La01, p.1; Kn98, p.59; An76, p.1].
A {\it partition}~$\l$ is a sequence of positive 
integers $\l=(\l_1, \l_2,\cdots, \l_\ell)$ such that 
$\l_1\geq \l_2 \geq \cdots \geq \l_\ell>0$.
The integers
$(\l_i)_{i=1,2,\ldots, \ell}$ are called the {\it parts} of~$\l$,
the number $\ell$ of parts being the
{\it length} of $\l$ denoted by $\ell(\l)$.  
The sum of its parts $\l_1+ \l_2+\cdots+ \l_\ell$ is
denoted by $|\l|$.
Let $n$ be an integer, a partition 
$\l$ is said to be a partition of $n$ if $|\l|=n$. We write $\l\vdash n$.
The set of all partitions of $n$
is denoted by $\setP(n)$. 
The set of all partitions is denoted by~$\setP$,
so that $$\setP=\bigcup_{n\geq 0} \setP(n).$$
Each partition can be represented by its Ferrers diagram. For example,
$\l=(6,3,3,2)$ is a partition  and its Ferrers diagram is reproduced in 
Fig.~\sec.1.

{
\long\def\maplebegin#1\mapleend{}

\maplebegin

# --------------- begin maple ----------------------

# Copy the following text  to "makefig.mpl"
# then in maple > read("makefig.mpl");
# it will create a file "z_fig_by_maple.tex"

#\unitlength=1pt

Hu:= 12.4; # height quantities
Lu:= Hu; # large unity

X0:=-95.0; Y0:=15.6; # origin position

File:=fopen("z_fig_by_maple.tex", WRITE);

mhook:=proc(x,y,lenx, leny)
local i, d,sp, yy, xx, ct;
	sp:=Hu/8;
	ct:=0;
	for xx from x*Lu+X0 to x*Lu+X0+Hu by sp do
		yy := y*Hu+Y0; 
		fprintf(File, "\\vline(
				xx,   yy+sp*ct-0.2,    Lu*leny-sp*ct+0.1);
		ct:=ct+1;
	od:
	
	ct:=0;
	for yy from y*Hu+Y0 to y*Hu+Y0+Hu by sp do
		xx := x*Lu+X0; 
		fprintf(File, "\\hline(
				xx+sp*ct-0.2,   yy,    Lu*lenx-sp*ct+0.1);
		ct:=ct+1;
	od:

end;

mhook(1,1,2,3);

fclose(File);
# -------------------- end maple -------------------------
\mapleend

\setbox1=\hbox{$
\def\b{\\{\hbox{}}}
\ytableau{12pt}{0.4pt}{}
{\b &\b       \cr 
 \b &\b &\b   \cr
 \b &\b &\b   \cr
 \b &\b &\b &\b &\b &\b  \cr
\noalign{\vskip 3pt}
\noalign{\hbox{Fig.~\sec.1. Partition}}
}$
}
\setbox2=\hbox{$
\def\b{\\{\hbox{}}}
\unitlength=1pt%
\def\vline(#1,#2)#3|{\leftput(#1,#2){\lline(0,1){#3}}}%
\def\hline(#1,#2)#3|{\leftput(#1,#2){\lline(1,0){#3}}}%
\ytableau{12pt}{0.4pt}{}
{\b &\b       \cr 
 \b &\b &\b   \cr
 \b &\b &\b   \cr
 \b &\b &\b &\b &\b &\b  \cr
\noalign{\vskip 3pt}
\noalign{\hbox{Fig. \sec.2. Hook length}}%
}
\vline(-82.6,27.8)37.3|
\vline(-81.0,29.4)35.8|
\vline(-79.5,30.9)34.2|
\vline(-78.0,32.4)32.6|
\vline(-76.4,34.0)31.1|
\vline(-74.8,35.6)29.6|
\vline(-73.3,37.1)28.0|
\vline(-71.8,38.6)26.4|
\vline(-70.2,40.2)24.9|
\hline(-82.8,28.0)24.9|
\hline(-81.2,29.6)23.4|
\hline(-79.7,31.1)21.8|
\hline(-78.2,32.6)20.2|
\hline(-76.6,34.2)18.7|
\hline(-75.0,35.8)17.2|
\hline(-73.5,37.3)15.6|
\hline(-72.0,38.8)14.0|
\hline(-70.4,40.4)12.5|
$
}
\setbox3=\hbox{$
\ytableau{12pt}{0.4pt}{}
{\\2 &\\1       \cr 
 \\4 &\\3 &\\1   \cr
 \\5 &\\4 &\\2   \cr
 \\9 &\\8 &\\6 &\\3 &\\2 &\\1  \cr
\noalign{\vskip 3pt}
\noalign{\hbox{Fig. \sec.3. Hook lengths}}
}$
}
$$\box1\quad\box2\quad\box3$$
}

For each box $v$ in the Ferrers diagram of a partition $\l$, or
for each box $v$ in $\l$, for short, define the 
{\it hook length} of $v$, denoted by $h_v(\l)$ or $h_v$, to be the number of 
boxes $u$ such that  $u=v$,
or $u$ lies in the same column as $v$ and above $v$, or in the 
same row as $v$ and to the right of $v$ (see Fig.~\sec.2). 
The {\it hook length multi-set} of $\l$, denoted by $\setH(\l)$,
is the multi-set of all hook lengths of $\l$.
Let $t$ be a positive integer. 
We write
$$
{\setH}_t(\l)=\{h\mid h\in\setH(\l), h\equiv 0 (\mod t)\}.
$$
In Fig.~\sec.3 
the hook lengths of all boxes for the partition $\l=(6,3,3,2)$
have been written in each box.
We have $\setH(\l)=\{2,1,4,3,1,5,4,2,9,8,6,\discretionary{}{}{} 3, 2, 1\}$
and
$\setH_2(\l)=\{2,4,4,2,8,6,2\}$.
\medskip

Recall that a partition $\l$ is a $t$-core if the {\it hook length multi-set} 
of $\l$ does not contain the integer $t$.
It is known that the hook length multi-set of each $t$-core 
does not contain any {\it multiple} of $t$ 
[Kn98. p.69, p.612; St99, p.468; JK81, p.75]. In other words, a partiton
$\l$ is a $t$-core if and only if $\setH_t(\l)=\emptyset$.

\medskip

{\it Definition \sec.1}.
Let $t=2t'+1$ be an odd positive integer. 
Each vector of integers $(v_0, v_1, \ldots, v_{t-1})\in \setZ^t$ is called 
{\it $V$-coding} if the following conditions hold:
(i) $v_i\equiv i (\mod t)$ for $0\leq i\leq t-1$; 
(ii) $v_0+v_1+\cdots + v_{t-1}=0$.

\smallskip
The $V$-coding is implicitly introduced in [Ma72]. It can be identified with 
the {\it set} 
$\{ v_0, v_1,  \ldots, v_{t-1}\}$ thanks to condition (i).
\medskip
Our first result is the following property on $t$-cores, which can be seen
as a refinement of the Nekrasov-Okounkov formula. 
The proof of this property involves 
a bijection 
between $t$-cores and integer vectors constructed by Garvan, Kim and 
Stanton [GKS90].

\proclaim Theorem \sec.1.
Let $t=2t'+1$ be an odd positive integer. 
There is a bijection $\phi_V: \l\mapsto (v_0, v_1, \ldots, v_{t-1})$ 
which maps each $t$-core onto a $V$-coding such that
$$
|\l|={1\over 2t}(v_0^2+v_1^2+\cdots +v_{t-1}^2) - {t^2-1\over 24}
\leqno{(\sec.1)}
$$
and
$$
\prod_{v\in\l} \Bigl(1-{t^2\over h_v^2}\Bigr)
={(-1)^{t'}\over 1!\cdot 2!\cdot 3!\cdots (t-1)!} 
\prod_{0\leq i<j\leq t-1} (v_i-v_j).
\leqno{(\sec.2)}
$$

We will describe the bijection $\phi_V$ and
prove the two equalities (\sec.1) and (\sec.2) in Section 2.
An example is given after
the construction of the bijection $\phi_V$.

\medskip
Next we provide an elementary proof of the following hook length formula,
discovered by
Nekrasov and Okounkov in their study of the Seiberg-Witten Theory 
[NO06, formula (6.12)].  
Our proof is based on the Macdonald identities for $A_\ell^{(a)}$ [Ma72]
and Theorem \sec.1.

\proclaim Theorem \sec.2 [Nekrasov-Okounkov].
For any complex number $z$ we have 
$$
\sum_{\l\in \setP}x^{|\l|} \prod_{h\in\setH(\l)}\bigl(1-{z \over h^2}\bigr)
\ =\ 
\prod_{k\geq 1} { (1-x^k)^{z-1}}.
\leqno{(\sec.3)}
$$

\medskip

Then we prove the following $(t,y)$-extension of Theorem \sec.2. When $y=t=1$
in (\sec.4) we recover the Nekrasov-Okounkov formula.
This extension unifies the Macdonald identities and the generating function for
$t$-cores.

\proclaim Theorem \sec.3. 
Let $t$ be a positive integer. For any complex numbers $y$ and $z$ 
we have 
$$
\sum_{\l\in\setP} x^{|\l|} \prod_{h\in \setH_t(\l)} 
\bigl(y-{tyz\over h^2}\bigr) 
= 
\prod_{k\geq 1}
{
(1-x^{tk})^t
\over
(1-(yx^t)^k)^{t-z}(1-x^k)  
}.\leqno{(\sec.4)}
$$

The proof of Theorem \sec.3, given in Section 4, is based 
on the 
Nekrasov-Okounkov formula (\sec.3) 
and on the properties of a classical bijection
which maps each partition to its $t$-core and $t$-quotient
[Ma95, p.12; St99, p.468; JK81, p.75; GSK90]. The following
result has a similar proof. 

\proclaim Theorem \sec.4. 
For any complex number $y$ we have 
$$
\sum_{\l\in\setP} x^{|\l|} y^{\#\{h\in\setH(\l), h=t\}} 
= 
\prod_{k\geq 1}
{
(1+(y-1)x^{tk})^t
\over
1-x^k
}.\leqno{(\sec.5)}
$$

\medskip
Last, we derive several applications of Theorems \sec.2 and \sec.3.
Let us single out some of them in this introduction. 
See [Ha08a] (resp. Section~5) for other applications of Theorem \sec.2 
(resp. Theorem \sec.3).

\proclaim Theorem \sec.5 [marked hook formula].
We have
$$
\sum_{\l\vdash n} f_\l^2\sum_{h\in\setH(\l)} h^2 = {n(3n-1)\over 2} n!,
\leqno{(\sec.6)}
$$
where $f_\l$ is the number of standard Young tableaux of shape $\l$.

The two sides of (\sec.6) can be combinatorially interpreted (see [Ha08a]).
However, a natural bijection between those two sides remains to be constructed.
Theorem \sec.5 is to be compared with the following well-known formula,
which is also a consequence of the
Robinson-Schensted-Knuth 
correspondence (see, for example, [Kn98, p.49-59; St99, p.324]).

$$
\sum_{\l\vdash n} f_\l^2 =  n! 
\leqno{(\sec.7)}
$$

The following theorem, proved in Section 6, 
improves a result due to Kostant [Ko04].

\proclaim Theorem \sec.6.
Let $k$ be a positive integer and $s$ be a real number 
such that $s\geq k^2-1$.  Then $(-1)^k f_k(s)>0$,  where
$f_k(s)$ is defined by
$$
\prod_{n\geq 1} (1-x^n)^s  =  \sum_{k\geq 0} f_k(s) x^k. 
$$

In section 7 we study the reversion of the Euler Product and obtain,
in particular, the following result. 

\proclaim Theorem \sec.7.
For any positive integers $n$ and $k$  the following two expressions
$$
\sum_{\l\vdash n} \prod_{v\in\l} \bigl(1+{k\over h_v^2}\bigr)
\leqno{(\sec.8)}
$$
and
$$
{1\over n+1}\sum_{\l\vdash n} \prod_{v\in\l} \bigl(1+{n\over h_v^2}\bigr)
\leqno{(\sec.9)}
$$
are integers.

The following specializations have 
similar forms, namely, Corollaries \sec.8, \sec.9 and \sec.10 on the one hand,
Corollaries \sec.11 and \sec.12 on the other hand. 
In fact, our motivation for Theorem \sec.3 was to look for a formula that
could interpolate the following two formulas (\sec.10) and (\sec.11).

\proclaim Corollary \sec.8 [\hbox{\rm $y=t=1$, $z=t^2$ in Theorem \sec.3}].
We have 
$$
\sum_{\l}x^{|\l|} \prod_{h\in\setH(\l)}\bigl(1-{t^2 \over h^2}\bigr)
\ =\ 
\prod_{k\geq 1} { (1-x^k)^{t^2} \over 1-x^k},
\leqno{(\sec.10)}
$$
where the sum ranges over all $t$-cores. 
\goodbreak

\proclaim Corollary \sec.9 [\hbox{\rm $z=t$ {\rm or } $y=0$ in 
Theorem \sec.3}]. 
We have 
$$
\sum_\l x^{|\l|}
=\prod_{k\geq 1} {(1-x^{tk} )^t\over 1-x^k}, \leqno{(\sec.11)}
$$
where the sum ranges over all $t$-cores. 
\goodbreak

Note that identity (\sec.11) is the well-known generating function
for $t$-cores [Ma95, p.12; St99, p.468; GSK90]. It is also the special case
$y=0$ of Theorem \sec.4.
The following 
identity is similar to the above two identities.
It is also 
a consequence of Theorem \sec.4.

\proclaim Corollary \sec.10 [\hbox{\rm $y=2$ in Theorem \sec.4}]. 
We have 
$$
\sum_{\l\in\setP} x^{|\l|} 2^{\#\{h\in\setH(\l), h=t\}} 
= 
\prod_{k\geq 1}
{
(1+x^{tk})^t
\over
1-x^k
}.\leqno{(\sec.12)}
$$

\proclaim Corollary \sec.11 [\hbox{\rm $y=t=1, z=2$ in Theorem \sec.3}]. 
We have 
$$
\sum_{\l\in\setP} x^{|\l|} \prod_{h\in \setH(\l)} 
\bigl(1-{2\over h^2}   \bigr) 
= 
\prod_{k\geq 1} (1-x^k). \leqno{(\sec.13)}
$$

\proclaim Corollary \sec.12 [\hbox{\rm $t=2, y=z=1$ in Theorem \sec.3}]. 
We have 
$$
\sum_{\l\in\setP} x^{|\l|} \prod_{h\in \setH_2(\l)} 
\bigl(1-{2\over h^2}   \bigr) 
= 
\prod_{k\geq 1} (1+x^k). \leqno{(\sec.14)}
$$

We end the introduction with some remarks.
The right-hand side of (\sec.13) can be expanded
by using the Euler pentagonal theorem [Eu83; An76, p.11]
$$
\prod_{k\geq 1}  (1-x^k)
= \sum_{m=-\infty}^{\infty}(-1)^m x^{m(3m+1)/2},
\leqno{(\sec.15)}
$$
so that Corollary \sec.11 says that
$$
\sum_{\l\vdash n} \prod_{h\in \setH(\l)} 
\bigl(1-{2\over h^2}   \bigr)  \leqno{(\sec.16)}
$$
is equal to $-1, 0, 1$ depending on the numerical value of $n$. 
\medskip

The right-hand side of (\sec.14) is the generating
function for partitions with {\it distinct} 
parts, so that Corollary \sec.12 says that
$$
\sum_{\l\vdash n} \prod_{h\in \setH_2(\l)} 
\bigl(1-{2\over h^2}   \bigr) \leqno{(\sec.17)}
$$
is equal to the number of partitions of $n$ with distinct parts.
\medskip

For example, there are five partitons of $n=4$ and two of them have 
distinct parts.
{
\medskip
\setbox1=\hbox{$\ytableau{12pt}{0.4pt}{}
{\\1   \cr 
 \\2   \cr
 \\3   \cr
 \\4   \cr
}$}
\setbox2=\hbox{$\ytableau{12pt}{0.4pt}{}
{\\1    \cr 
 \\2    \cr
 \\4 &\\1  \cr
}$}
\setbox3=\hbox{$\ytableau{12pt}{0.4pt}{}
{\\2 &\\1    \cr
 \\3 &\\2  \cr
}$}
\setbox4=\hbox{$\ytableau{12pt}{0.4pt}{}
{\\1    \cr 
 \\4 &\\2 &\\1  \cr
}$}
\setbox5=\hbox{$\ytableau{12pt}{0.4pt}{}
{ \\4 &\\3 &\\2 &\\1  \cr
}$}
\centerline{\box1\qquad\box2\qquad\box3\qquad\box4\qquad\box5}
\nobreak \smallskip \nobreak
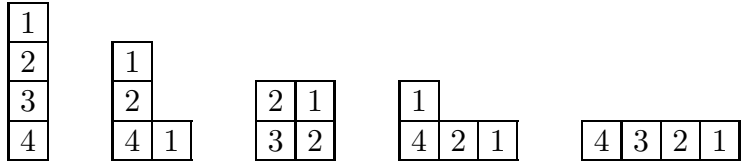
\centerline{Fig. \sec.4. The multi-set of hook lengths for $\setP(4)$}
\medskip
}
\noindent
We have
{
\def\fh#1|{\bigl(1-{2\over #1^2}\bigr)}
$$
\leqalignno{
  & 2\fh1|\fh2|\fh3|\fh4|\cr
+ & 2\fh1|\fh1|\fh2|\fh4|\cr
+ & \fh1|\fh2|\fh2|\fh3|=0\cr
}
$$
}
and
{
\def\fh#1|{\bigl(1-{2\over #1^2}\bigr)}
$$ 2\fh2|\fh4| +  2\fh2|\fh4| +  \fh2|\fh2|=2. $$
}

It would be interesting to explain directly why (\sec.16) and (\sec.17)
are integers.
\def\sec{2}
\section{\sec. New properties of $t$-cores} 

In this section 
we first describe the bijection $\phi_V$ required in Theorem 1.1
and then prove equalities (1.1) and (1.2).
Let $t=2t'+1$ be an odd positive integer. 
Each finite set of integers $A=\{a_1, a_2, \ldots, a_n\}$ is said to be 
{\it $t$-compact} if the following conditions hold:

(i) $-1, -2, \ldots, -t\in A$;

(ii) for each $a\in A$ such that 
$a\not=-1,-2, \ldots, -t$, we have $a\geq 1$ and $a \not\equiv 0 \mod t$; 

(iii) let $b>a\geq 1$ be two integers such that $a \equiv b \mod t$. 
If $b\in A$,  then $a\in A$.
\smallskip

Let $A$ be a $t$-compact set. An element $a\in A$ is said to be 
{\it $t$-maximal}
if $b\not\in A$ for every $b>a$ such that $a\equiv b\mod t$.
The set of $t$-maximal elements of $A$ is denoted by
$\tmax(A)$. 
Let $\l$ be a $t$-core.
The {\it $H$-set} of the $t$-core $\l$ is defined to be 
$$
H(\l)=\{h_v \mid \hbox{$v$ is a box in the leftmost 
column of $\l$}\}\cup\{-1,-2, \ldots -t\}.
$$

The notion of $H$-set is a variation of the {\it $\beta$-numbers} introduced
by James and Kerber, who also introduced
the runners-beads-abacus model [JK81, p.75] in the study of $t$-cores.
In this section we prefer to work directly on the $H$-sets, as our goal is
to prove identities (1.1) and (1.2).

\medskip

\proclaim Lemma \sec.1.
For each $t$-core $\l$ its $H$-set 
$H(\l)$ is a $t$-compact set.

{\it Proof}. Let $c=tk+r$ ($k\geq 1, 0\leq r\leq t-1$) be an element in $H(\l)$ 
and $a$ be the maximal element in $H(\l)$ such that $a<t(k-1)+r$. 
We must show that $t(k-1)+r$ is also in $H(\l)$. 
If it were not the case, let 
$z>t(k-1)+r, y_1, y_2, \ldots, y_d$ be the hook lengths as shown in 
Fig. \sec.1, 
where only the relevant horizontal section of the 
partition diagram has been represented.
We have
$y_1=c-a-1 \geq tk+r - t(k-1)-r  = t$
and 
$y_d=c-z+1 \leq tk+r - t(k-1)-r = t$;
so that there is one hook $y_i=t$. 
This is a contradiction since $\l$
is supposed to be a $t$-core. 
\qed

\medskip 
{
\setbox1=\hbox{$
\def\b{\\{\hbox{}}}
\ytableau{14pt}{0.4pt}{}
{
 \\a  &\b  &\b  &\b  &\b \cr
 \\z   &\b  &\b  &\b  &\b  &\b  &\b &\b  &\b  \cr
 \b   &\b  &\b  &\b  &\b  &\b  &\b &\b  &\b  &\b \cr
 \b   &\b  &\b  &\b  &\b  &\b  &\b &\b  &\b  &\b \cr
 \\c  &\b  &\b  &\b  &\b &\\{y_1}  &\\{y_2} &\\{\cdots}  &\\{y_d} 
     &\b &\b  &\b  \cr
\noalign{\vskip 3pt}
\noalign{\hbox{Fig. \sec.1. Hook length and $t$-compact set}}
}$
}
\midinsert
$$
\box1
$$
\endinsert
}

{\it Construction of $\phi_V$}. Let $\l$ be a $t$-core and $H(\l)$ be
its $H$-set.
The {\it $U$-coding} of $\l$ is defined to be the set
$U:=\tmax(H(\l))$, which can be identified with the vector
$(u_0, u_1, \ldots, u_{t-1})$ such that 
$u_0=-t$, $u_i>-t$ and $u_i\equiv i \mod t$ for $1\leq i\leq t-1$.
In general,
$$S:=u_0 + u_1 +\cdots + u_{t-1} \not=0.\leqno{(\sec.1)}$$ 
The integer $S$ is a multiple of $t$ because
$$
S=\sum u_i = \sum (tk_i +i) = t\sum k_i + t(t-1)/2 \leqno{(\sec.2)}
$$
(remember that $t=2t'+1$ is an odd integer). 
The $V$-coding $\phi_V(\l)$ is the set $V$ obtained from $U$ by 
the following {\it normalization}: 
$$
\phi_V(\l)=V:=\{u-S/t \ : \ u\in U\}. \leqno{(\sec.3)}
$$
In fact, we can prove that $S/t=\ell(\l)-t'-1$  (see (\sec.8)).
The set $V$ can be identified with a vector $V$-coding because 
$$\sum v_i=\sum (u_i -S/t)=\sum u_i - S=0.$$
\medskip

\noindent
{\it Example \sec.1}.
Consider the $5$-core 
$$\l=(14,10,6,6,4,4,4,2,2,2).$$ 
The $H$-set of $\l$ (see Fig. \sec.2)
$$
H(\l)=\{23,18,13,12,9,8,7,4,3,2, -1, -2, -3, -4, -5\}
$$
is $5$-compact. 
The $U$-coding of $\l$ is $U=\max_5(H(\l))=\{23,12,9,-4,-5\}$,
or in vector form  
$$(u_0,u_1,u_2,u_3,u_4)=(-5, -4, 12, 23, 9).$$
As $S=\sum u_i=35$,  the $V$-coding is given by
$$V=\{ -5-7, -4-7, 12-7, 23-7, 9-7 \}
=\{ -12, -11, 5, 16, 2 \},$$
or in vector form
$$
\phi_V(\l)=(v_0, v_1, v_2, v_3, v_4)=(5, 16, 2, -12, -11).
$$



\long\def\maplebegin#1\mapleend{}

\maplebegin

# --------------- begin maple ----------------------

# Copy the following text  to "makefig.mpl"
# then in maple > read("makefig.mpl");
# it will create a file "z_fig_by_maple.tex"

#\unitlength=1pt

Hu:= 16.4; # height quantities
Lu:= Hu; # large unity

X0:=0; Y0:=40; # origin position

File:=fopen("z_fig_by_maple.tex", WRITE);

hdash:=proc(x,y,len)
local i, d,sp, xx;
	d:=2;
	sp:=2;
	for xx from x*Lu+X0 to x*Lu+X0+Lu*len-d by d+sp do
	fprintf(File, "\\hline(
	od:
end;

vdash:=proc(x,y,len)
local i, d,sp, yy;
	d:=2;
	sp:=2;
	for yy from y*Hu+Y0 to y*Hu+Y0+Hu*len-d by d+sp do
		fprintf(File, "\\vline(
	od:
end;

text:=proc(x,y,t)
	fprintf(File, "\\centerput(
end;

cercle:=proc(x,y) text(x,y, "\\cerclechar"); end;
# the partition

# h-line

cercle(0,14);    cercle(-1,14);   
cercle(0,13);    cercle(-1,13);
cercle(0,4);     cercle(3,4);
cercle(0,3);     cercle(5,3);
cercle(0,0);     cercle(13,0);

text(-1,14, "{\\it 0}");text(0,14, "-5");
text(-1,13, "{\\it 1}");text(0,13, "-4");
text(-1,12, "{\\it 2}");text(0,12, "-3");
text(-1,11, "{\\it 3}");text(0,11, "-2");
text(-1,10, "{\\it 4}");text(0,10, "-1");
text(-1,9, "{\\it 0}");

hdash(-1,0,1);
hdash(-1,15,2);
vdash(-1,0,15);
vdash(0,10,5);
vdash(1,10,5);

# region separ
hhdash:=proc(x,y) hdash(x,y+0.03,1); end;
HHdash:=proc(XL, dx) local x; for x in XL do hhdash(x, x+dx); od; end;
vvdash:=proc(x,y) vdash(x,y+0.03,1); end;
VVdash:=proc(XL, dx) local x; for x in XL do vvdash(x, x+dx); od; end;

HHdash([-1,0,1,2,3], 11);
VVdash([   0,1,2,3], 10);

HHdash([-1,0,1,2,3,4,5,6], 6);
VVdash([   0,1,2,3,4,5,6], 5);

HHdash([-1,0,1,2,3,4,5,6,7,8], 1);
VVdash([   0,1,2,3,4,5,6,7,8], 0);

HHdash([4,5,6,7,8,9,10,11], -4);
VVdash([  5,6,7,8,9,10,11], -5);

HHdash([9,10,11,12,13], -9);
VVdash([  10,11,12,13], -10);

text(2,14.2, "$r=-2$");
text(4,12.2, "$r=-1$");
text(6,9.8, "$r=0$");
text(8,6.8, "$r=1$");
text(11,3.8, "$r=2$");
text(12.5,1.5, "$r=3$");

fclose(File);
# -------------------- end maple -------------------------
\mapleend


{
\unitlength=1pt
\font\cerclefont=cmsy10 at 14pt
\def\cerclechar{\hbox{\cerclefont\char'015}}
\def\vline(#1,#2)#3|{\leftput(#1,#2){\lline(0,1){#3}}}
\def\hline(#1,#2)#3|{\leftput(#1,#2){\lline(1,0){#3}}}
\newbox\boxhook
\setbox\boxhook=\vbox{
\vskip 26mm
$
\def\b{\\{\hbox{}}}
\def\r#1{\\{{\it #1}}}
\hskip 48pt \ytableau{16pt}{0.4pt}{}
{
 \\2 &\r2       \cr 
 \\3 &\r3       \cr 
 \\4 &\r4       \cr 
 \\7 &\r0 &\r1 &\r2  \cr
 \\8 &\b &\b &\r3  \cr
 \\9 &\b &\b &\r4  \cr
 \\{12} &\b &\b &\r0 &\r1 &\r2 \cr
 \\{13} &\b &\b &\b &\b &\r3 \cr
 \\{18} &\b &\b &\b &\b &\r4 &\r0 &\r1 &\r2 &\r3 \cr
 \\{23} &\b &\b &\b &\b &\b &\b &\b
        &\b &\r4 &\r0 &\r1 &\r2 &\r3 \cr
\noalign{\vskip 3pt}
\noalign{\hbox{\quad Fig. \sec.2. $U$-coding and $N$-coding of $t$-core}}
\centerput(0,270){\cerclechar}
\centerput(-16,270){\cerclechar}
\centerput(0,254){\cerclechar}
\centerput(-16,254){\cerclechar}
\centerput(0,106){\cerclechar}
\centerput(50,106){\cerclechar}
\centerput(0,90){\cerclechar}
\centerput(82,90){\cerclechar}
\centerput(0,40){\cerclechar}
\centerput(214,40){\cerclechar}
\centerput(-16,270){{\it 0}}
\centerput(0,270){-5}
\centerput(-16,254){{\it 1}}
\centerput(0,254){-4}
\centerput(-16,237){{\it 2}}
\centerput(0,237){-3}
\centerput(-16,221){{\it 3}}
\centerput(0,221){-2}
\centerput(-16,204){{\it 4}}
\centerput(0,204){-1}
\centerput(-16,188){{\it 0}}
\hline(-24,36)2|
\hline(-20,36)2|
\hline(-16,36)2|
\hline(-12,36)2|
\hline(-24,282)2|
\hline(-20,282)2|
\hline(-16,282)2|
\hline(-12,282)2|
\hline(-8,282)2|
\hline(-4,282)2|
\hline(0,282)2|
\hline(4,282)2|
\vline(-24,36)2|
\vline(-24,40)2|
\vline(-24,44)2|
\vline(-24,48)2|
\vline(-24,52)2|
\vline(-24,56)2|
\vline(-24,60)2|
\vline(-24,64)2|
\vline(-24,68)2|
\vline(-24,72)2|
\vline(-24,76)2|
\vline(-24,80)2|
\vline(-24,84)2|
\vline(-24,88)2|
\vline(-24,92)2|
\vline(-24,96)2|
\vline(-24,100)2|
\vline(-24,104)2|
\vline(-24,108)2|
\vline(-24,112)2|
\vline(-24,116)2|
\vline(-24,120)2|
\vline(-24,124)2|
\vline(-24,128)2|
\vline(-24,132)2|
\vline(-24,136)2|
\vline(-24,140)2|
\vline(-24,144)2|
\vline(-24,148)2|
\vline(-24,152)2|
\vline(-24,156)2|
\vline(-24,160)2|
\vline(-24,164)2|
\vline(-24,168)2|
\vline(-24,172)2|
\vline(-24,176)2|
\vline(-24,180)2|
\vline(-24,184)2|
\vline(-24,188)2|
\vline(-24,192)2|
\vline(-24,196)2|
\vline(-24,200)2|
\vline(-24,204)2|
\vline(-24,208)2|
\vline(-24,212)2|
\vline(-24,216)2|
\vline(-24,220)2|
\vline(-24,224)2|
\vline(-24,228)2|
\vline(-24,232)2|
\vline(-24,236)2|
\vline(-24,240)2|
\vline(-24,244)2|
\vline(-24,248)2|
\vline(-24,252)2|
\vline(-24,256)2|
\vline(-24,260)2|
\vline(-24,264)2|
\vline(-24,268)2|
\vline(-24,272)2|
\vline(-24,276)2|
\vline(-24,280)2|
\vline(-8,200)2|
\vline(-8,204)2|
\vline(-8,208)2|
\vline(-8,212)2|
\vline(-8,216)2|
\vline(-8,220)2|
\vline(-8,224)2|
\vline(-8,228)2|
\vline(-8,232)2|
\vline(-8,236)2|
\vline(-8,240)2|
\vline(-8,244)2|
\vline(-8,248)2|
\vline(-8,252)2|
\vline(-8,256)2|
\vline(-8,260)2|
\vline(-8,264)2|
\vline(-8,268)2|
\vline(-8,272)2|
\vline(-8,276)2|
\vline(-8,280)2|
\vline(9,200)2|
\vline(9,204)2|
\vline(9,208)2|
\vline(9,212)2|
\vline(9,216)2|
\vline(9,220)2|
\vline(9,224)2|
\vline(9,228)2|
\vline(9,232)2|
\vline(9,236)2|
\vline(9,240)2|
\vline(9,244)2|
\vline(9,248)2|
\vline(9,252)2|
\vline(9,256)2|
\vline(9,260)2|
\vline(9,264)2|
\vline(9,268)2|
\vline(9,272)2|
\vline(9,276)2|
\vline(9,280)2|
\hline(-24,201)2|
\hline(-20,201)2|
\hline(-16,201)2|
\hline(-12,201)2|
\hline(-8,217)2|
\hline(-4,217)2|
\hline(0,217)2|
\hline(4,217)2|
\hline(9,234)2|
\hline(13,234)2|
\hline(17,234)2|
\hline(21,234)2|
\hline(25,250)2|
\hline(29,250)2|
\hline(33,250)2|
\hline(37,250)2|
\hline(42,267)2|
\hline(46,267)2|
\hline(50,267)2|
\hline(54,267)2|
\vline(-8,201)2|
\vline(-8,205)2|
\vline(-8,209)2|
\vline(-8,213)2|
\vline(9,217)2|
\vline(9,221)2|
\vline(9,225)2|
\vline(9,229)2|
\vline(25,234)2|
\vline(25,238)2|
\vline(25,242)2|
\vline(25,246)2|
\vline(42,250)2|
\vline(42,254)2|
\vline(42,258)2|
\vline(42,262)2|
\hline(-24,119)2|
\hline(-20,119)2|
\hline(-16,119)2|
\hline(-12,119)2|
\hline(-8,135)2|
\hline(-4,135)2|
\hline(0,135)2|
\hline(4,135)2|
\hline(9,152)2|
\hline(13,152)2|
\hline(17,152)2|
\hline(21,152)2|
\hline(25,168)2|
\hline(29,168)2|
\hline(33,168)2|
\hline(37,168)2|
\hline(42,185)2|
\hline(46,185)2|
\hline(50,185)2|
\hline(54,185)2|
\hline(58,201)2|
\hline(62,201)2|
\hline(66,201)2|
\hline(70,201)2|
\hline(74,217)2|
\hline(78,217)2|
\hline(82,217)2|
\hline(86,217)2|
\hline(91,234)2|
\hline(95,234)2|
\hline(99,234)2|
\hline(103,234)2|
\vline(-8,119)2|
\vline(-8,123)2|
\vline(-8,127)2|
\vline(-8,131)2|
\vline(9,135)2|
\vline(9,139)2|
\vline(9,143)2|
\vline(9,147)2|
\vline(25,152)2|
\vline(25,156)2|
\vline(25,160)2|
\vline(25,164)2|
\vline(42,168)2|
\vline(42,172)2|
\vline(42,176)2|
\vline(42,180)2|
\vline(58,185)2|
\vline(58,189)2|
\vline(58,193)2|
\vline(58,197)2|
\vline(74,201)2|
\vline(74,205)2|
\vline(74,209)2|
\vline(74,213)2|
\vline(91,217)2|
\vline(91,221)2|
\vline(91,225)2|
\vline(91,229)2|
\hline(-24,37)2|
\hline(-20,37)2|
\hline(-16,37)2|
\hline(-12,37)2|
\hline(-8,53)2|
\hline(-4,53)2|
\hline(0,53)2|
\hline(4,53)2|
\hline(9,70)2|
\hline(13,70)2|
\hline(17,70)2|
\hline(21,70)2|
\hline(25,86)2|
\hline(29,86)2|
\hline(33,86)2|
\hline(37,86)2|
\hline(42,103)2|
\hline(46,103)2|
\hline(50,103)2|
\hline(54,103)2|
\hline(58,119)2|
\hline(62,119)2|
\hline(66,119)2|
\hline(70,119)2|
\hline(74,135)2|
\hline(78,135)2|
\hline(82,135)2|
\hline(86,135)2|
\hline(91,152)2|
\hline(95,152)2|
\hline(99,152)2|
\hline(103,152)2|
\hline(107,168)2|
\hline(111,168)2|
\hline(115,168)2|
\hline(119,168)2|
\hline(124,185)2|
\hline(128,185)2|
\hline(132,185)2|
\hline(136,185)2|
\vline(-8,37)2|
\vline(-8,41)2|
\vline(-8,45)2|
\vline(-8,49)2|
\vline(9,53)2|
\vline(9,57)2|
\vline(9,61)2|
\vline(9,65)2|
\vline(25,70)2|
\vline(25,74)2|
\vline(25,78)2|
\vline(25,82)2|
\vline(42,86)2|
\vline(42,90)2|
\vline(42,94)2|
\vline(42,98)2|
\vline(58,103)2|
\vline(58,107)2|
\vline(58,111)2|
\vline(58,115)2|
\vline(74,119)2|
\vline(74,123)2|
\vline(74,127)2|
\vline(74,131)2|
\vline(91,135)2|
\vline(91,139)2|
\vline(91,143)2|
\vline(91,147)2|
\vline(107,152)2|
\vline(107,156)2|
\vline(107,160)2|
\vline(107,164)2|
\vline(124,168)2|
\vline(124,172)2|
\vline(124,176)2|
\vline(124,180)2|
\hline(58,37)2|
\hline(62,37)2|
\hline(66,37)2|
\hline(70,37)2|
\hline(74,53)2|
\hline(78,53)2|
\hline(82,53)2|
\hline(86,53)2|
\hline(91,70)2|
\hline(95,70)2|
\hline(99,70)2|
\hline(103,70)2|
\hline(107,86)2|
\hline(111,86)2|
\hline(115,86)2|
\hline(119,86)2|
\hline(124,103)2|
\hline(128,103)2|
\hline(132,103)2|
\hline(136,103)2|
\hline(140,119)2|
\hline(144,119)2|
\hline(148,119)2|
\hline(152,119)2|
\hline(156,135)2|
\hline(160,135)2|
\hline(164,135)2|
\hline(168,135)2|
\hline(173,152)2|
\hline(177,152)2|
\hline(181,152)2|
\hline(185,152)2|
\vline(74,37)2|
\vline(74,41)2|
\vline(74,45)2|
\vline(74,49)2|
\vline(91,53)2|
\vline(91,57)2|
\vline(91,61)2|
\vline(91,65)2|
\vline(107,70)2|
\vline(107,74)2|
\vline(107,78)2|
\vline(107,82)2|
\vline(124,86)2|
\vline(124,90)2|
\vline(124,94)2|
\vline(124,98)2|
\vline(140,103)2|
\vline(140,107)2|
\vline(140,111)2|
\vline(140,115)2|
\vline(156,119)2|
\vline(156,123)2|
\vline(156,127)2|
\vline(156,131)2|
\vline(173,135)2|
\vline(173,139)2|
\vline(173,143)2|
\vline(173,147)2|
\hline(140,37)2|
\hline(144,37)2|
\hline(148,37)2|
\hline(152,37)2|
\hline(156,53)2|
\hline(160,53)2|
\hline(164,53)2|
\hline(168,53)2|
\hline(173,70)2|
\hline(177,70)2|
\hline(181,70)2|
\hline(185,70)2|
\hline(189,86)2|
\hline(193,86)2|
\hline(197,86)2|
\hline(201,86)2|
\hline(206,103)2|
\hline(210,103)2|
\hline(214,103)2|
\hline(218,103)2|
\vline(156,37)2|
\vline(156,41)2|
\vline(156,45)2|
\vline(156,49)2|
\vline(173,53)2|
\vline(173,57)2|
\vline(173,61)2|
\vline(173,65)2|
\vline(189,70)2|
\vline(189,74)2|
\vline(189,78)2|
\vline(189,82)2|
\vline(206,86)2|
\vline(206,90)2|
\vline(206,94)2|
\vline(206,98)2|
\centerput(33,273){$r=-2$}
\centerput(66,241){$r=-1$}
\centerput(99,201){$r=0$}
\centerput(132,152){$r=1$}
\centerput(181,103){$r=2$}
\centerput(205,65){$r=3$}
}$
}
\midinsert
$$
\box\boxhook
$$
\vskip -20pt
\endinsert
}
\noindent
We have
$$
\leqalignno{
|\l|&={1\over 2t}(v_0^2+v_1^2+\cdots +v_{t-1}^2) - {t^2-1\over 24}\cr
&= {1\over 2\cdot 5} ( 5^2+ 16^2+ 2^2+(-12)^2+(-11)^2) - {5^2-1\over 24}=54.\cr
}
$$
and
$$
\leqalignno{
\prod_{v\in\l} \Bigl(1-{5^2\over h_v^2}\Bigr)
&={1\over 1!\cdot 2!\cdot 3!\cdots (t-1)!} 
\prod_{0\leq i<j\leq t-1} (v_i-v_j)\cr
&=(-11)(3)(17)(16)\cdot(14)(28)(27)\cdot(14)(13)\cdot(-1)/288\cr
&=60035976.\cr
}
$$
Notice that, as expected, the above two numbers are positive integers.
\medskip
A vector of integers $(n_0, n_1, \ldots, n_{t-1})\in \setZ^t$ is said to be an 
{\it $N$-coding} if $n_0+n_1+\cdots+n_{t-1}=0$. 
Garvan, Kim and Stanton have defined a bijection $\phi_N$ between
$N$-codings and $t$-cores. We now recall its definition using their own words 
[GKS90,p.3] (see also [BG06]).

Let $\l$ be a $t$-core. Define the vector $(n_0, \ldots, n_{t-1})=\phi_N(\l)$
in the following way. Label the box in the $i$-th row and $j$-column of $\l$
by $j-i \mod t$. 
We also label the boxes in column 0 (in dotted lines in Fig. \sec.2) in
the same way, and call the resulting diagram the {\it extended $t$-residue
diagram}. A box is called {\it exposed} if it is at the end of a row of the
extended $t$-residue diagram. 
The set of boxes $(i,j)$ satisfying 
$t(r-1)\leq j-i < tr$ of the extended $t$-residue diagram of $\l$ is 
called {\it region} and numbered $r$. 
In Fig. \sec.2 the regions have been bordered by 
dotted lines.  
We now define $n_i$ to be the 
maximum region $r$ which contains an exposed box labeled $i$.

In Fig. \sec.2 the labels of all boxes lying on the maximal border strip
(but the leftmost one) have been written in italic. This includes all the
exposed boxes: 3,3,3,2,4,3,2,4,3,2,4,3,2,1,0, when reading from bottom to top.
We have 
$(n_0, n_1, n_2, n_3, n_4)=(-2, -2, 1, 3, 0)$.

\proclaim Theorem \sec.2 [Garvan-Kim-Stanton].
The bijection 
$$\phi_N : \l \mapsto (n_0, n_1, \ldots, n_{t-1})$$
has the following property:
$$
|\l| = {t\over 2} \sum_{i=0}^{t-1} n_i^2 + \sum_{i=0}^{t-1} in_i. 
\leqno{(\sec.4)}
$$

Let $t'=(t-1)/2$ and let
$$\phi_V^N : (n_0, n_1, \ldots, n_{t-1}) \mapsto
(v_0, v_1, \ldots, v_{t-1})$$
be the bijection that maps each $N$-coding onto the $V$-coding defined by
$$
v_i=\cases
{
t n_{i+t'} +i & if $0\leq i\leq t'$; \cr
t n_{i-t'-1} +i-t & if $t'+1\leq i\leq t-1$ \cr
}
\leqno{(\sec.5)}
$$
or in set form
$$
\{v_i \mid 0\leq i \leq t-1\}= \{ t n_i +i -t' \mid 0\leq i \leq t-1\}.
\leqno{(\sec.6)}
$$
The bijective property of the map $\phi_V^N$ is easy to verify.
More essentially, the bijection $\phi_V$ is the composition 
product of the two previous bijections as is now shown.

\proclaim Lemma \sec.3.
We have $\phi_V = \phi_V^N\circ\phi_N$.

{\it Proof}. 
Let $(v_0, \ldots, v_{t-1})=\phi_V(\l)$,
$(n_0, \ldots, n_{t-1})=\phi_N(\l)$ and
$$(v_0', \ldots, v_{t-1}')=\phi_V^N(n_0, \ldots, n_{t-1}).$$
We need prove that $v_i=v_i'$.
The number $n_i$ in the $N$-coding is
defined to be the maximum region $r$ which contains an exposed box 
labelled $i$. This exposed box is called {\it critical italic box}. 
In Fig. \sec.2 a circle is drawn around the label 
of each critical italic box.
On the other hand,
the $U$-coding is defined to be the set 
$\tmax(H(\l))$, where $H(\l)$ is the $H$-set of~$\l$.
A box in the leftmost column whose hook length is an element of the $U$-coding
is called {\it critical roman box}. In Fig. \sec.2, a circle is drawn
around the hook length number of each critical roman box.
Let us write the labels of all the exposed boxes (the vector $L=(L_i)$) with
its region numbers (the vector $R=(R_i)$)
and the $H$-set of $\l$ (the vector $H=(H_i)=H(\l)$), 
read from bottom to top.
{$$%
\font\cerclefont=cmsy10 at 14pt%
\def\a#1{#1\kern-9pt\hbox{\cerclefont\char'015}}%
\def\b#1{#1\kern-12pt\hbox{\cerclefont\char'015}}%
\def\c#1{#1\kern-11pt\hbox{\cerclefont\char'015}}%
\def\ng#1{\!\hbox{\rm -#1}}%
\matrix{
L=&\a{\it3} &\it 3 &\it 3 &\a{\it2}   &\a{\it4}&\it3&\it2&\it4&\it3&
\it 2& \it4& \it3& \it2& \a{\it1}&\a{\it0} & \cr
R=&\a{3}&2&1&\a{1}&\a0&0&0&\ng1&\ng1&\ng1&\ng2&\ng2&\ng2&\c{\ng2}&\c{\ng2}& \cr
H=&\b{23}&18&13&\b{12}&\a9&8&7&4&3&2&\ng1&\ng2&\ng3&\c{\ng4}&\c{\ng5}& \cr
}%
$$}%
It is easy to see that $L_j= (H_j-\ell(\l))\mod t$ 
and
$R_j= \lfloor(H_j-\ell(\l))/t\rfloor +1$.
This means that $L_i$ has a circle symbol if and only if $H_i$ has a 
circle symbol.
We then have a natural bijection
$$f: u_i \mapsto  \lfloor(u_i-\ell(\l))/t\rfloor +1 = n_{(u_i-\ell) \mod t}
\leqno{(\sec.7)}
$$
between the set $\{u_0, \ldots, u_{t-1}\}$ and 
$\{n_0, \ldots, n_{t-1}\}$. 
By (\sec.6) and (\sec.7) we have
$$
\leqalignno{
\{v_i'\} 
& = \{ tn_i+i-t'\} \cr
& = \{ t n_{ (u_i-\ell) \mod t} + (u_i-\ell) \mod t - t'\} \cr
& = \{ t( \lfloor(u_i-\ell)/t\rfloor +1  ) + (u_i-\ell) \mod t - t'\} \cr
& = \{ u_i-\ell+t'+1 \}. \cr
}
$$
On the other hand, $(v_i')$ is a $V$-coding,
because $v_i'\equiv i\mod t$ and 
$\sum v_i' = t\sum n_i + \sum i -t (t-1)/2 = 0$;
so that 
$$(\sum_i u_i )/t = \ell -t'-1. \leqno{(\sec.8)}$$
Hence
$$
\{v_i'\} 
 = \{ u_i-\ell+t'+1 \}
 = \{ u_i-(\sum_i u_i)/t \} 
 = \{ v_i \}.\qed 
$$

\medskip
Take again the same partition as in Example \sec.1; the $N$-coding
is
$$(n_0, n_1, n_2, n_3, n_4)=(-2, -2, 1, 3, 0).$$
We verify that
$$
\leqalignno{
&(v_0', v_1', v_2', v_3', v_4')\cr
&\qquad=(1\times5+0,\ 3\times5+1,\ 0\times5+2,\ -2\times5-2,\ -2\times5-1).\cr
&\qquad=(5, 16, 2, -12, -11)=(v_0, v_1, v_2, v_3, v_4).\cr
}
$$
\medskip

{\it Proof of (1.1) in Theorem 1.1}.
From (\sec.6) we have
$$
\leqalignno{
\sum v_i^2&=\sum(tn_i+i-t')^2\cr
&=\sum \bigl((tn_i)^2 +2tin_i -2tt'n_i +i^2 +t'^2 -2it'  \bigr) \cr
&= t^2\sum n_i^2 + 2t\sum in_i + {(t-1)t(2t-1)\over 6} + tt'^2 - t' t (t-1)\cr
&= t^2\sum n_i^2 + 2t\sum in_i  + {t(t^2-1)\over 12}.
}
$$
Hence
$$
{1\over 2t}\sum v_i^2 = {t\over 2}\sum  n_i^2 + \sum in_i + {t^2-1\over 24} 
= |\l| + {t^2-1\over 24}.\qed
$$
\medskip
For proving (1.2) in Theorem 1.1, we first etablish the following two lemmas.
\proclaim Lemma \sec.4. 
For any $t$-compact set $A$ we have
$$
\prod_{a\in A, a>0} \Bigl(1-{t^2\over a^2} \Bigr)
=
\prod_{a\in \tmax(A), a\not=-t} {a+t\over a} .\leqno{(\sec.9)}
$$

\medskip

\noindent
{\it Example \sec.2}.
Take $t=5$. Then the set 
$$A=\{-5,-4, -3, -2, -1, 2,3,4,7,8,9,12,13,18,23\}$$ is $5$-compact.
We have $\tmax(A)=\{-5,-4, 9,12,23\}$. Hence
$$
\prod_{a\in A, a>0} \Bigl(1-{25\over a^2} \Bigr)=
{ 1\cdot 14 \cdot 17\cdot 28\over (-4) \cdot 9\cdot 12\cdot 23 }.
\leqno{(\sec.10)}
$$
{\it Proof}. Write
\vskip -10pt
$$
\prod_{a\in A, a>0} \Bigl(1-{t^2\over a^2} \Bigr)
=
\prod_{a\in A, a>0} {(a-t)\cdot (a+t)\over a \cdot a},
$$
then
delete the common factors in numerator and denominator, as illustrated 
by means of Example \sec.2. 
$$
{\def\coef#1|#2|#3|{{{#1 \over #3} \times {#2\over  #3}}}%
\def\ccf#1|#2|{{{#1 \over #2} }}%
\matrix{
\ccf1|-4| & & & & & & (a\equiv 1\mod 5)\cr
\noalign{\medskip}
\ccf2|-3|  & \coef-3|7|2| &\coef2|12|7| &\coef7|17|12| && &(a\equiv 2\mod 5)\cr
\noalign{\medskip}
\ccf3|-2| &\coef-2|8|3| &\coef3|13|8| &\coef8|18|13| &\coef13|23|18| 
   &\coef18|28|23| &(a\equiv 3\mod 5)\cr
\noalign{\medskip}
  \ccf4|-1|& \coef-1|9|4|  &\coef4|14|9| &&& &(a\equiv 4\mod 5)\cr
}}
$$
The product $(a-5)(a+5)/a^2$ for $a>0$ is reproduced in the row
determined by $a\mod 5$ in the above table, except for the leftmost column.
But the product of the factors in the leftmost column is equal to $1$ because
$t$ is an odd integer;
so that
the left-hand side of (\sec.10) is the product of the factors in the above
table. After deleting the common factors, it remains the rightmost 
fraction in each row.\qed

\medskip

\proclaim Lemma \sec.5.
Let $\l$ be a $t$-core and $(u_0, u_1, \ldots, u_{t-1})$ be its $U$-coding
(defined in the body of the construction of $\phi_V$).
Let $\l'$ be the $t$-core obtained from $\l$ by erasing the leftmost column
of $\l$ and $(u'_0, u'_1, \ldots, u'_{t-1})$ be its $U$-coding.
Then
$$
\prod_{0\leq i<j\leq t-1} {u_i-u_j \over u_i'-u_j'} = 
\prod_{j=1}^{t-1} {u_j+t\over u_j}. 
$$

\medskip
\noindent
{\it Example \sec.3}. Take the 5-core $\l$ given in Example \sec.1. The
$U$-coding of $\l$ is $(u_0,u_1,u_2,u_3,u_4)=(-5, -4, 12, 23, 9)$.
We have 
$$\l'=(13,9,5,5,3,3,3,1,1,1).$$ 
The $U$-coding of $\l'$ is
$(u'_0,u'_1,u'_2,u'_3,u'_4)=(-5, 11, 22, 8, -1)$.
Now, consider the cyclic rearrangement
$$(u''_0,u''_1,u''_2,u''_3,u''_4)=(-1, -5, 11, 22, 8)$$
 of $(u'_0,u'_1,u'_2,u'_3,u'_4)$.
We have $\prod(u'_i-u'_j)=\prod(u''_i-u''_j)$ because $t$ is an odd integer.
Moreover $u''_i=u_i-1$ for all $1\leq i\leq 4$. 
Hence
$$
\eqalignno{
 \prod_{0\leq i<j\leq t-1} { u_i-u_j \over u''_i-u''_j}
&=
 \prod_{j=1}^{t-1} { u_0-u_j \over u''_0-u''_j}\cr
&= { (-5+4)(-5-12)(-5-23)(-5-9) \over (-1+5)(-1-11)(-1-22)(-1-8) }\cr
&= { (-4+5)(12+5)(23+5)(9+5) \over (-4)(12)(23)(9) }.\cr
}
$$

{\it Proof}.
We suppose that $\l$ contains $\delta$ parts equal to $1$.
Its $H$-set  $H(\l)$ 
(viewed as a vector in decreasing order if necessary) 
can be split into six segments
$H(\l)=A_1A_2A_3A_4A_5A_6$ defined by (see Fig. \sec.3)

(i) $a \geq \delta+2$ for each $a\in A_1$;

(ii) $A_2=(\delta, \delta-1, \ldots, 3,2,1)$;

(iii) $A_3=(-1, -2, -3, \ldots, \delta+2-t)$;

(iv) $A_4=(\delta+1-t)$;

(v) $A_5=(\delta-t, \delta-1-t, \ldots, 1-t)$;

(vi) $A_6=(-t)$.

\smallskip

On the other hand the $H$-set $H(\l')$ of $\l'$
is split into five segments
$H(\l')=A_1'A_2'A_3'A_4'A_5'$ defined by 

(i') $A_1'=\{a-\delta-1 \ : \ a\in A_1\}$;

(ii') $A_2'=\{a-\delta-1 \ : \ a\in A_2\}=(-1, -2, \ldots, -\delta)$;

(iii') $A_3'=(-\delta-1)$;

(iv') $A_4'=\{a-\delta-1 \ : \ a\in A_3\}=(-\delta-2, -\delta-3, \ldots, -t+1)$;

(v) $A_5'=(-t)$.

%


\long\def\maplebegin#1\mapleend{}

\maplebegin

# --------------- begin maple ----------------------

# Copy the following text  to "makefig.mpl"
# then in maple > read("makefig.mpl");
# it will create a file "z_fig_by_maple.tex"

Hu:= 4;  # height unity
Hn:= 17; # height quantities
Lu:= 20; # large unity

X0:=-34; Y0:=0; # origin position

File:=fopen("z_fig_by_maple.tex", WRITE);

hline:=proc(x,y,len)
	fprintf(File, "\\hline(
end;

vline:=proc(x,y,len)
	fprintf(File, "\\vline(
end;

hdash:=proc(x,y,len)
local i, d,sp, xx;
	d:=1;
	sp:=1;
	for xx from x*Lu+X0 to x*Lu+X0+Lu*len-d by d+sp do
	fprintf(File, "\\hline(
	od:
end;

vdash:=proc(x,y,len)
local i, d,sp, yy;
	d:=1;
	sp:=1;
	for yy from y*Hu+Y0 to y*Hu+Y0+Hu*len-d by d+sp do
		fprintf(File, "\\vline(
	od:
end;

text:=proc(x,y,t)
	fprintf(File, "\\centerput(
			y*Hu+Y0+1-Hu,    t);
end;

textup:=proc(x,y,t)
	fprintf(File, "\\centerput(
			y*Hu+Y0+1-ceil(Hu/2),    t);
end;

# the partition

# h-line

hdash(0,18,1);
hdash(0,17,1); 
               hdash(1,14,1);
hdash(0,13,1); hdash(1,13,1);
hdash(0,12,1);
               hdash(1,9,1);
hline(0,8,1);  hdash(1,8,1);
hdash(0,4,1);  hline(1,4,1);
hdash(2,0,2);  hline(2,2,2);
hline(0,0,4);

# v-line

vdash(0,8,10);
                 vdash(1,8,10);
                 vline(1,4,4);  vdash(2,4,10);
                                vline(2,2,2);
vline(0,0,8);    vdash(1,0,4);  vdash(2,0,2);      vline(4,0,2);

# text
text(1,18, "$-t$");
text(1,17, "$1-t$");
text(1,16, "$\\vdots$");
text(1,15, "$\\delta-1-t$");
text(1,14, "$\\delta-t$");    text(2,14, "$-t$");
text(1,13, "$\\delta+1-t$");  text(2,13, "$1-t$");
text(1,12, "$\\delta+2-t$");  text(2,12, "$\\vdots$");
text(1,11, "$\\vdots$");      text(2,11, "$-\\delta-3$");
text(1,10, "$-2$");           text(2,10, "$-\\delta-2$");
text(1,9, "$-1$");            text(2,9, "$-\\delta-1$");
text(1,8, "$1$");             text(2,8, "$-\\delta$");
text(1,7, "$2$");             text(2,7, "$\\vdots$");
text(1,6, "$\\vdots$");       text(2,6, "$-2$");
text(1,5, "$\\delta$");       text(2,5, "$-1$");                  
#
text(0, 18, "$A_6$");
textup(0, 15, "$A_5$");  text(3, 14, "$A_5'$");
text(0, 13, "$A_4$");    textup(3, 11, "$A_4'$");
textup(0, 10, "$A_3$");  text(3, 9, "$A_3'$");
textup(0, 6, "$A_2$");   textup(3, 6, "$A_2'$");
textup(1, 2, "$A_1$");   textup(2, 2, "$A_1'$"); 
#                    

fclose(File);

# -------------------- end maple -------------------------
\mapleend


\newbox\boxarbre
\def\vline(#1,#2)#3|{\leftput(#1,#2){\lline(0,1){#3}}}
\def\hline(#1,#2)#3|{\leftput(#1,#2){\lline(1,0){#3}}}
\setbox\boxarbre=\vbox{\vskip
70mm\offinterlineskip 
%
\hline(-34,72)1|
\hline(-32,72)1|
\hline(-30,72)1|
\hline(-28,72)1|
\hline(-26,72)1|
\hline(-24,72)1|
\hline(-22,72)1|
\hline(-20,72)1|
\hline(-18,72)1|
\hline(-16,72)1|
\hline(-34,68)1|
\hline(-32,68)1|
\hline(-30,68)1|
\hline(-28,68)1|
\hline(-26,68)1|
\hline(-24,68)1|
\hline(-22,68)1|
\hline(-20,68)1|
\hline(-18,68)1|
\hline(-16,68)1|
\hline(-14,56)1|
\hline(-12,56)1|
\hline(-10,56)1|
\hline(-8,56)1|
\hline(-6,56)1|
\hline(-4,56)1|
\hline(-2,56)1|
\hline(0,56)1|
\hline(2,56)1|
\hline(4,56)1|
\hline(-34,52)1|
\hline(-32,52)1|
\hline(-30,52)1|
\hline(-28,52)1|
\hline(-26,52)1|
\hline(-24,52)1|
\hline(-22,52)1|
\hline(-20,52)1|
\hline(-18,52)1|
\hline(-16,52)1|
\hline(-14,52)1|
\hline(-12,52)1|
\hline(-10,52)1|
\hline(-8,52)1|
\hline(-6,52)1|
\hline(-4,52)1|
\hline(-2,52)1|
\hline(0,52)1|
\hline(2,52)1|
\hline(4,52)1|
\hline(-34,48)1|
\hline(-32,48)1|
\hline(-30,48)1|
\hline(-28,48)1|
\hline(-26,48)1|
\hline(-24,48)1|
\hline(-22,48)1|
\hline(-20,48)1|
\hline(-18,48)1|
\hline(-16,48)1|
\hline(-14,36)1|
\hline(-12,36)1|
\hline(-10,36)1|
\hline(-8,36)1|
\hline(-6,36)1|
\hline(-4,36)1|
\hline(-2,36)1|
\hline(0,36)1|
\hline(2,36)1|
\hline(4,36)1|
\hline(-34,32)20|
\hline(-14,32)1|
\hline(-12,32)1|
\hline(-10,32)1|
\hline(-8,32)1|
\hline(-6,32)1|
\hline(-4,32)1|
\hline(-2,32)1|
\hline(0,32)1|
\hline(2,32)1|
\hline(4,32)1|
\hline(-34,16)1|
\hline(-32,16)1|
\hline(-30,16)1|
\hline(-28,16)1|
\hline(-26,16)1|
\hline(-24,16)1|
\hline(-22,16)1|
\hline(-20,16)1|
\hline(-18,16)1|
\hline(-16,16)1|
\hline(-14,16)20|
\hline(6,0)1|
\hline(8,0)1|
\hline(10,0)1|
\hline(12,0)1|
\hline(14,0)1|
\hline(16,0)1|
\hline(18,0)1|
\hline(20,0)1|
\hline(22,0)1|
\hline(24,0)1|
\hline(26,0)1|
\hline(28,0)1|
\hline(30,0)1|
\hline(32,0)1|
\hline(34,0)1|
\hline(36,0)1|
\hline(38,0)1|
\hline(40,0)1|
\hline(42,0)1|
\hline(44,0)1|
\hline(6,8)40|
\hline(-34,0)80|
\vline(-34,32)1|
\vline(-34,34)1|
\vline(-34,36)1|
\vline(-34,38)1|
\vline(-34,40)1|
\vline(-34,42)1|
\vline(-34,44)1|
\vline(-34,46)1|
\vline(-34,48)1|
\vline(-34,50)1|
\vline(-34,52)1|
\vline(-34,54)1|
\vline(-34,56)1|
\vline(-34,58)1|
\vline(-34,60)1|
\vline(-34,62)1|
\vline(-34,64)1|
\vline(-34,66)1|
\vline(-34,68)1|
\vline(-34,70)1|
\vline(-14,32)1|
\vline(-14,34)1|
\vline(-14,36)1|
\vline(-14,38)1|
\vline(-14,40)1|
\vline(-14,42)1|
\vline(-14,44)1|
\vline(-14,46)1|
\vline(-14,48)1|
\vline(-14,50)1|
\vline(-14,52)1|
\vline(-14,54)1|
\vline(-14,56)1|
\vline(-14,58)1|
\vline(-14,60)1|
\vline(-14,62)1|
\vline(-14,64)1|
\vline(-14,66)1|
\vline(-14,68)1|
\vline(-14,70)1|
\vline(-14,16)16|
\vline(6,16)1|
\vline(6,18)1|
\vline(6,20)1|
\vline(6,22)1|
\vline(6,24)1|
\vline(6,26)1|
\vline(6,28)1|
\vline(6,30)1|
\vline(6,32)1|
\vline(6,34)1|
\vline(6,36)1|
\vline(6,38)1|
\vline(6,40)1|
\vline(6,42)1|
\vline(6,44)1|
\vline(6,46)1|
\vline(6,48)1|
\vline(6,50)1|
\vline(6,52)1|
\vline(6,54)1|
\vline(6,8)8|
\vline(-34,0)32|
\vline(-14,0)1|
\vline(-14,2)1|
\vline(-14,4)1|
\vline(-14,6)1|
\vline(-14,8)1|
\vline(-14,10)1|
\vline(-14,12)1|
\vline(-14,14)1|
\vline(6,0)1|
\vline(6,2)1|
\vline(6,4)1|
\vline(6,6)1|
\vline(46,0)8|
\centerput(-24,69){$-t$}
\centerput(-24,65){$1-t$}
\centerput(-24,61){$\vdots$}
\centerput(-24,57){$\delta-1-t$}
\centerput(-24,53){$\delta-t$}
\centerput(-4,53){$-t$}
\centerput(-24,49){$\delta+1-t$}
\centerput(-4,49){$1-t$}
\centerput(-24,45){$\delta+2-t$}
\centerput(-4,45){$\vdots$}
\centerput(-24,41){$\vdots$}
\centerput(-4,41){$-\delta-3$}
\centerput(-24,37){$-2$}
\centerput(-4,37){$-\delta-2$}
\centerput(-24,33){$-1$}
\centerput(-4,33){$-\delta-1$}
\centerput(-24,29){$1$}
\centerput(-4,29){$-\delta$}
\centerput(-24,25){$2$}
\centerput(-4,25){$\vdots$}
\centerput(-24,21){$\vdots$}
\centerput(-4,21){$-2$}
\centerput(-24,17){$\delta$}
\centerput(-4,17){$-1$}
\centerput(-44,69){$A_6$}
\centerput(-44,59){$A_5$}
\centerput(16,53){$A_5'$}
\centerput(-44,49){$A_4$}
\centerput(16,43){$A_4'$}
\centerput(-44,39){$A_3$}
\centerput(16,33){$A_3'$}
\centerput(-44,23){$A_2$}
\centerput(16,23){$A_2'$}
\centerput(-24,7){$A_1$}
\centerput(-4,7){$A_1'$}
%
%
}
\midinsert
$$
\box\boxarbre
$$
\centerline{\qquad\quad Fig. \sec.3. 
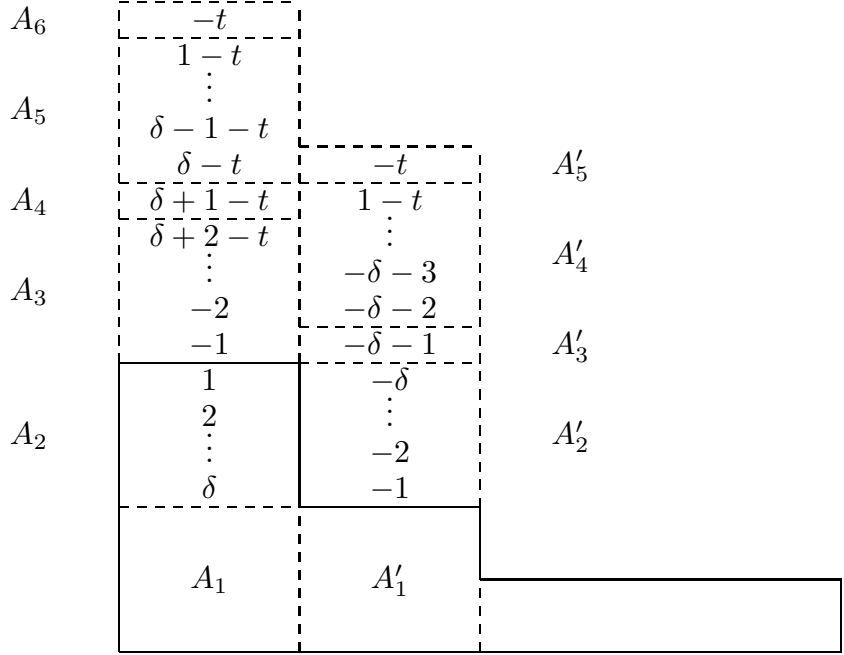
Comparison of the hook lengths of $\l$ and $\l'$}
\endinsert
\smallskip
\noindent
Notice that some segments $A_i$ and $A_i'$ may be empty. More
precisely, 
$$
\cases{
A_2=A_5=A_2'=\emptyset, &  if $\delta=0$; \cr
A_3=A_4'=\emptyset, &  if $\delta=t-2$; \cr
A_3=A_4=A_3'=A_4'=\emptyset, & if $\delta=t-1$. \cr
}
$$

The basic facts are:

(i) $a\not \in \tmax(H(\l))$ for every $a\in A_5$; 
because 
$\{a\mod t\ : \ a\in A_5\} =\{a\mod t\ : \ a\in A_3\}$. 
In other words the set $A_5$  is {\it masked} by $A_3$. 

(ii) $\delta+1-t\in\tmax(H(\l))$; because $a\not\equiv 0\mod t$ for every 
$a\in A_1'$ so that
$a\not\equiv \delta+1 \mod t$ for every $a\in A_1$.
It is easy to see that
$a\not\equiv \delta+1 \mod t$ for every $a\in A_2\cup A_3$.

(iii) $-\delta-1\in\tmax(H(\l'))$; because $a\not\equiv 0\mod t$ for every
$a\in A_1 \cup A_2$ so that
$a\not\equiv -\delta-1\mod t$ for every $a\in A_1'\cup A_2'$.

(iv) Since that $a\mapsto a-\delta-1$ is a bijection between
$A_1\cup A_2\cup A_3$ and $A'_1\cup A'_2\cup A'_4$, 
it is also a bijection between
$\tmax(H(\l))\setminus\{-t,\delta-t+1\}$
and $\tmax(H(\l'))\setminus\{-t,-\delta-1\}$.
\medskip

The above facts enable us to derive the $U$-coding of $\l'$ from the
$U$-coding of $\l$ as follows. Let
$$(u_i)=(u_0=-t, u_1, u_2, \ldots, u_{k-1},  \delta+1-t, 
u_{k+1}, u_{k+1}, \ldots, u_{t-1})$$
be the $U$-coding of $\l$ and define
$$
(u''_i)=(u''_0=-\delta-1, u_1'', u_2'', \ldots, u_{k-1}'',  -t, 
u_{k+1}'', u_{k+1}'', \ldots, u_{t-1}'')$$
where $u''_i=u_i-\delta-1$ for $i\geq 1$.
Then, the $U$-coding of $\l'$ is simply
$$
(u'_i)=(u'_0=-t,  u_{k+1}'', u_{k+1}'', \ldots, u_{t-1}'',
-\delta-1, u_1'', u_2'', \ldots, u_{k-1}'').
$$
We have $\prod(u'_i-u'_j)=\prod(u''_i-u''_j)$ because $t$ is an odd integer.
On the other hand, 
$u''_i-u''_j=u_i-u_j$ for all $1\leq i<j\leq t-1$. Hence 
$$
\eqalignno{
 \prod_{0\leq i<j\leq t-1} { u_i-u_j \over u'_i-u'_j}
&= \prod_{0\leq i<j\leq t-1} { u_i-u_j \over u''_i-u''_j}
= \prod_{j=1}^{t-1} { u_0-u_j \over u''_0-u''_j}\cr
&= \prod_{j=1}^{t-1} { -t-u_j \over -\delta-1-u''_j}
= \prod_{j=1}^{t-1} { u_j +t\over u_j}.\qed\cr
}
$$
\medskip

{\it Proof of (1.2) in Theorem 1.1}.
Because the $U$-coding and $V$-coding of~$\l$ only differ by the
normalization given in (\sec.3) and $t$ is an odd integer, we have
$\prod(v_i-v_j)=\prod(u_i-u_j)$.
By Lemmas \sec.5 and \sec.4 we have 
$$
\leqalignno{
\prod_{0\leq i<j\leq t-1} (u_i-u_j) 
&=\prod_{j=1}^{t-1} {u_j+t\over u_j}\times 
  \prod_{0\leq i<j\leq t-1} (u'_i-u'_j) \cr
&= 
\prod_{a\in H(\l), a>0} \Bigl(1-{t^2\over a^2} \Bigr)\times
  \prod_{0\leq i<j\leq t-1} (u'_i-u'_j) \cr
&= \cdots 
=
K\times \prod_{v\in \l} \Bigl(1-{t^2\over h_v^2} \Bigr). \cr
}
$$
Taking $\l$ as the empty $t$-core, the $U$-coding of $\l$ is
$(-t, -t+1, -t+2, \ldots, -3, -2, -1)$. 
We then obtain $K= (-1)^{t'}{ 1!\cdot 2!\cdot 3!\cdots (t-1)!} $\qed
\def\sec{3}
\section{\sec. Expansion formula for the powers of the Euler Product} 
The powers of the Euler Product and the hook lengths of partitions 
are two mathematical objects widely studied in the Theory of
Partitions, in Algebraic Combinatorics and Group Representation Theory.
In this section we give an elementary proof of Theorem 1.2, which 
establishes a new connection by giving an explicit 
expansion formula for all the powers $s$ of the Euler
Product in terms of partition hook lengths, where the exponent $s$
is any complex number. Recall that the {\it Euler Product}
is the infinite product $\prod_{m\geq 0} (1-x^m)$. 
The following two formulas [Eu83; An76, p.11, p.21]
go back to Euler (the pentagonal theorem)
$$
\prod_{m\geq 1}  (1-x^m)
= \sum_{k=-\infty}^{\infty}(-1)^k x^{k(3k+1)/2}
\leqno{(\sec.1)}
$$
and Jacobi (triple product identity, 
see [An76, p.21; Kn98, p.20; JS89; FH99; FK99])
$$
\prod_{m\geq 1} (1-x^m)^3=\sum_{m\geq 0} (-1)^m (2m+1)x^{m(m+1)/2}.
\leqno{(\sec.2)}
$$
Further explicit formulas for the
powers of the Euler Product
$$
\prod_{m\geq 1} (1-x^m)^s  =  \sum_{k\geq 0} f_k(s) x^k \leqno{(\sec.3})
$$
have been derived for certain integers
$$s=1,3,8,10,14,15,21,24,26,28,35,36,\ldots \leqno{(\sec.4)}$$
by 
Klein and Fricke for $s=8$, 
Atkin for $s=14, 26$,  Winquist for $s=10$, 
and Dyson for $s=24, \ldots$ [Wi69; Dy72].  
The paper entitled ``Affine root systems and Dedekind's $\eta$-function",
written by Macdonald in 1972, is a milestone in the study of powers
of Euler Product [Ma72].  The review of this paper for MathSciNet,
written by Verma [Ve], contains seven pages!  
It has also inspired several followers, see
[Ka74; Mo75; Ko76; Ko04; Mi85; AF02; CFP05; RS06].
The main achievement of
Macdonald was to unify all the 
well-known formulas for the integers $s$ listed in (\sec.4), 
except for $s=26$. 
He obtained an expansion formula of 
$$\prod_{m\geq 0} (1-x^m)^{\dim \lig}\leqno{(\sec.5)}$$
for every semi-simple Lie algebra $\lig$. 
\medskip
A variation of the Euler Product, called the {\it Dedekind $\eta$-function}, 
is defined by 
$$
\eta(x)=x^{1/24} \prod_{m\geq 1} (1-x^m). \leqno{(\sec.6)}
$$
We are ready to state the Macdonald identities for $A_\ell^{(a)}$ [Ma72],
which play a fundamental role in the following proof of Thorem 1.2.

\proclaim Theorem \sec.1 [Macdonald].
Let $t=2t'+1$ be an odd integer. 
We have
$$
\eta(x)^{t^2-1} = c_0 \sum_{(v_0, \ldots, v_{t-1})} \prod_{i<j} (v_i-v_j) 
x^{(v_0^2+v_1^2+\cdots+v_{t-1}^2)/(2t)},
\leqno{(\sec.7)}
$$
where the sum ranges over all $V$-codings $(v_0, v_1, \ldots, v_{t-1})$ 
(see Definition 1.1)
and $c_0$ is a numerical constant.

Consider the term of lowest degree in the above power series.
We immediately get 
$$
c_0={(-1)^{t'}\over 1!\cdot 2!\cdot 3!\cdots (t-1)!}. \leqno{(\sec.8)} 
$$
\medskip

{\it Proof of Theorem 1.2}.
Using the following identity (see [St99, p.316]) 
$$
\prod_{m\geq 1} {1\over 1-x^m} = 
\exp{\bigl(\sum_{k\geq 1} {x^k\over k(1-x^k)}\bigr)}, \leqno{(\sec.9)}
$$
the right-hand side of equation (1.3) can be written:
$$
\prod_{m\geq 1} {1\over 1-x^m}\times 
\exp{\Bigl(-z\sum_{k\geq 1} {x^k\over k(1-x^k)}\Bigr)}.
\leqno{(\sec.10)}
$$
Let $n\geq 0$ be a positive integer.
The coefficient $C_n(z)$ of $x^n$ on the left-hand side of 
(1.3) is a polynomial
in $z$ of degree $n$. The coefficient $D_n(z)$ of $x^n$ on the 
right-hand side
of (1.3) is also a 
polynomial in $z$ of degree $n$ thanks to (\sec.10). 
For  proving $C_n(z)=D_n(z)$, it suffices to find $n+1$
explicit numerical values $z_0, z_1, \ldots, z_n$ such that 
$C_n(z_i)= D_n(z_i)$ for $0\leq i\leq n$
by using the Lagrange interpolation formula.
The basic fact is that
$$
\prod_{v\in\l}\bigl(1-{t^2 \over h_v^2}\bigr)=0
$$
for every partition $\l$ which is not a $t$-core.
By comparing Theorems 1.1 and \sec.1 we see that equation (1.3)  is true
when $z=t^2$ for every odd integer~$t$, i.e.,
$$
\sum_{\l\in \setP}\ x^{|\l|} \prod_{v\in\l}\bigl(1-{t^2 \over h_v^2}\bigr)
\ =\ \prod_{m\geq 1} {(1-x^m)^{t^2-1}},
$$
so that $C_n(z)=D_n(z)$ for every complex number $z$.  \qed
\medskip

Note that Kostant already observed that $D_n(z)$ 
is a polynomial in~$z$, but did not mention any explicit expression [Ko04].

\def\sec{4}
\section{\sec. A unified hook formula via $t$-cores} 
In this section we prove Theorems 1.3 and 1.4 by using 
the properties of a classical bijection
which maps each partition to its $t$-core and $t$-quotient
[Ma95, p.12; St99, p.468; JK81, p.75; GSK90]. 
Let $\setW$ be the set of 
bi-infinite binary sequences beginning with infinitely
many 0's and ending with infinitely many 1's.
Each element $w$ of $\setW$ can be represented by 
$(b_i)_{i} =\cdots b_{-3}b_{-2}b_{-1}b_0b_1b_2b_3\cdots$,
but the representation is not unique. Actually, for any fixed integer $k$
the sequence $(b_{i+k})_i$ also represents $w$.
The {\it canonical representation} of $w$ is the unique sequence 
$(c_i)_i= \cdots c_{-3}c_{-2}c_{-1}c_0c_1c_2c_3\cdots$
such that 
$$
\#\{i\leq -1, c_i=1\}=\#\{i\geq 0, c_i=0\}.
$$
We put a dot symbol ``$.$" between the letters $c_{-1}$ and $c_0$ 
in the bi-infinite sequence $(c_i)_i$ when it is the 
canonical representation.
\medskip
There is a natural one-to-one correspondence between $\setP$ and $\setW$
(see, e.g. [St99, p.468; AF02] for more detail).
Let $\l$ be a partition. We encode each horizontal edge of $\l$ by 1 and each 
vertical edge by 0. Reading these (0,1)-encodings from top to bottom and from
left to right yields a binary word $u$. By adding infinitely many 0's 
to the left and infinitely many 1's to the right of $u$ we get
an element $w=\cdots 000u111\cdots\in\setW$. 
Clearly the map $\psi: \l\mapsto w$ is a one-to-one correspondence 
between $\setP$ and $\setW$.
The canonical representation of $\psi(\l)$ will be denoted by $C_\l$.
For example, take $\l=(6,5,3,3)$; we have $u=1110011010$, so that
$w=\cdots 0001110011010111 \cdots$ and 
$C_\l=(c_i)_i=\cdots 0001110.011010111 \cdots$
{
\long\def\maplebegin#1\mapleend{}%
\maplebegin

# --------------- begin maple ----------------------

# Copy the following text  to "makefig.mpl"
# then in maple > read("makefig.mpl");
# it will create a file "z_fig_by_maple.tex"

#\unitlength=1pt

Hu:= 16.4; # height quantities
Lu:= Hu; # large unity

X0:=-9.0*Lu; Y0:=Hu; # origin position

File:=fopen("z_fig_by_maple.tex", WRITE);

text:=proc(x,y,t)
	fprintf(File, "\\centerput(
end;

text0:=proc(x,y) text(x-0.35,y+0.6, "\\hbox{\\eightpoint 0}"); end;
text1:=proc(x,y) text(x,y+0.25, "\\hbox{\\eightpoint 1}"); end;

text1(1,4); text1(2,4); text1(3,4);
text1(4,2); text1(5,2);
text1(6,1);

text0(4,3);
text0(4,2);
text0(6,1);
text0(7,0);

fclose(File);

# -------------------- end maple -------------------------
\mapleend
\setbox2=\hbox{
$
\unitlength=1pt%
\ytableau{16pt}{0.4pt}{}
{
\noalign{\vskip 0pt}
 \\3 &\\2 &\\1  \cr
 \\4 &\\3 &\\2  \cr
 \\7 &\\6 &\\5 &\\2 &\\1   \cr
 \\9 &\\8 &\\7 &\\4 &\\3 &\\1   \cr
\noalign{\vskip 3pt}
\noalign{\hbox{\hskip -8.7mm Fig.~\sec.1. Partition and (0,1)-sequence}}
}
 \centerput(-131,87){\hbox{\eightpoint 1}}
 \centerput(-114,87){\hbox{\eightpoint 1}}
 \centerput(-98,87){\hbox{\eightpoint 1}}
 \centerput(-82,54){\hbox{\eightpoint 1}}
 \centerput(-65,54){\hbox{\eightpoint 1}}
 \centerput(-49,37){\hbox{\eightpoint 1}}
 \centerput(-87,76){\hbox{\eightpoint 0}}
 \centerput(-87,60){\hbox{\eightpoint 0}}
 \centerput(-54,43){\hbox{\eightpoint 0}}
 \centerput(-38,27){\hbox{\eightpoint 0}}
$}
\midinsert
$$\box2$$
\endinsert
}
%

Let $t$ be a positive integer. 
It is known [Ma95, p.12; St99, p.468; JK81, p.75; GSK90] that there is a 
bijection $\Omega$ which maps
a partition $\l$ to $(\mu; \l^0, \l^1, \ldots, \l^{t-1})$ such that
\smallskip 
{

(P1) $\mu$ is a $t$-core and $\l^0, \l^1, \ldots, \l^{t-1}$ are partitions;

(P2) $|\l|=|\mu|+t(|\l^0|+|\l^1|+\cdots+|\l^{t-1}|)$;

(P3) $\{h/t\mid h\in \setH_t(\l)\}=
\setH(\l^0)\cup \setH(\l^1)\cup \cdots \cup \setH(\l^{t-1})$.
}
\smallskip \noindent 
The vector $(\l^0, \l^1, \ldots, \l^{t-1})$ is 
usually called the {\it $t$-quotient}
of the partition~$\l$. Let us briefly describe the bijection $\Omega$ 
(see, e.g., [AF02; St99, p.468]).
We split the canonical representation $C_\l=(c_i)_i$ of the partition $\l$ 
into $t$ sections.  This means that we form the  subsequence
$w^{k}= (c_{it+k})_i$ for each $k=0,1,\ldots, t-1$.
The $k$-th entry $\l^k$ of the $t$-quotient of $\l$
is defined to be the inverse image $\psi^{-1}(w^k)$ of the subsequence $w^k$.
With the above example we have $w^0=\cdots 00110111\cdots$
and $w^1=00010100111\cdots$, so that $\l^0=(2,1)$ and $\l^1=(2,2,1)$.
Property (P3) holds since $\setH(\l^0)=\{2,1\}$, 
$\setH(\l^1)=\{1,3,1,4,2\}$ and
$\setH_2(\l)=\{2,4,2,6,2,8,4\}$ (See Fig. \sec.1-4).
Notice that the subsequence $w^k$ defined by
$w^{k}= (c_{it+k})_i$
is not necessarily the canonical representation. 
For that reason we do not reproduce
the dot symbol ``$.$" in the corresponding rows in the following tableau.

$$
\vbox{\offinterlineskip\halign{
\ \hfil$#$\hfil\ \vrule&&\strut\ \hfil$#$\hfil\ \cr
C_\l   & \cdots &0&0&0&0&0&1&1&1&0&.&0&1&1&0&1&0&1&1&1&\cdots \cr
\noalign{\hrule}
w^0   & \cdots  &&0& &0& &1& &1& & &0& &1& &1& &1& &1&\cdots \cr
v^0   & \cdots  &&0& &0& &0& &1& & &1& &1& &1& &1& &1&\cdots \cr
\noalign{\hrule}
w^1   & \cdots &0& &0& &0& &1& &0& & &1& &0& &0& &1& &\cdots \cr
v^1   & \cdots &0& &0& &0& &0& &0& & &0& &1& &1& &1& &\cdots \cr
\noalign{\hrule}
C_\mu  & \cdots &0&0&0&0&0&0&0&1&0&.&1&0&1&1&1&1&1&1&1&\cdots \cr
}}
$$
{
\setbox1=\hbox{$
\def\b{\\{\hbox{}}}
\ytableau{12pt}{0.4pt}{}
{\\2 &\\1       \cr 
\noalign{\vskip 3pt}
\noalign{\hbox{\hskip -3mm Fig.~\sec.2. Partition $\l^0$}}
}$
}
\setbox2=\hbox{$
\def\b{\\{\hbox{}}}
\unitlength=1pt%
\def\vline(#1,#2)#3|{\leftput(#1,#2){\lline(0,1){#3}}}%
\def\hline(#1,#2)#3|{\leftput(#1,#2){\lline(1,0){#3}}}%
\ytableau{12pt}{0.4pt}{}
{\\1        \cr 
 \\3 &\\1    \cr
 \\4 &\\2    \cr
\noalign{\vskip 3pt}
\noalign{\hbox{\hskip -3mm Fig.~\sec.3. Partition $\l^1$}}%
}
$
}
\setbox3=\hbox{$
\ytableau{12pt}{0.4pt}{}
{\\1       \cr 
 \\3 & \\1    \cr
\noalign{\vskip 3pt}
\noalign{\hbox{\hskip -3mm Fig.~\sec.4. The $2$-core $\mu$}}
}$
}
$$\quad\box1\qquad\box2\qquad\box3$$
}

For each subsequence $w^k$ we continually replace the subword $10$ by $01$.
The final resulting sequence is of the
form $\cdots 000111\cdots$ and is denoted by $v^k$.
The $t$-core  of the partition $\l$ is defined to be the partition $\mu$
such that the $t$ sections of the canonical representation $C_\mu$ 
are exactly $v^0, v^1, \ldots, v^{t-1}$. For the above example we have
$\mu=(2,1)$. Properties (P2) and (P3) can be derived from the following
basic fact: each box of $\l$ is in one-to-one correspondence with the
ordered pair of integers $(i,j)$ such that $i<j$ and $c_i=1, c_j=0$.
Moreover the hook length of that box is equal to $j-i$. 

\medskip

{\it Proof of Theorem 1.3}.
By the properties of the bijection $\Omega$ we get 
$$
\leqalignno{
\sum_{\l\in\setP} x^{|\l|} \prod_{h\in \setH_t(\l)} 
\bigl(y-{tyz\over h^2}\bigr) 
&=  \prod_{k\geq 1}{(1-x^{tk})^t\over 1-x^k}
\Bigl(
\sum_{\l\in\setP} x^{t|\l|} \prod_{h\in \setH(\l)} 
\bigl(y-{tyz\over (th)^2}\bigr) 
\Bigr)^t \cr
&=  \prod_{k\geq 1}{(1-x^{tk})^t\over 1-x^k}
\Bigl(
\sum_{\l\in\setP} (yx^t)^{|\l|} \prod_{h\in \setH(\l)} 
\bigl(1-{z/t\over h^2}\bigr) 
\Bigr)^t. &{(\sec.1)} \cr
}
$$
By Theorem 1.2 
$$\sum_{\l\in\setP} (yx^t)^{|\l|} \prod_{h\in \setH(\l)} 
\bigl(1-{z/t\over h^2}\bigr) 
=
\prod_{m\geq 1} { (1-(yx^t)^m)^{z/t-1}}. \leqno{(\sec.2)}
$$
We obtain (1.4) when reporting (\sec.2) into (\sec.1).\qed
\medskip

{\it Proof of Theorem 1.4}.
It is easy to see that
$$
\sum_{\l\in\setP} x^{|\l|} y^{\#\{h\in\setH(\l), h=1\}} 
= 
\prod_{m\geq 1}{1+(y-1)x^{m}\over 1-x^{m}}.
$$
By the properties of the bijection $\Omega$ we get 
$$
\leqalignno{
\sum_{\l\in\setP} x^{|\l|} y^{\#\{h\in\setH(\l), h=t\}} 
&=  \prod_{k\geq 1}{(1-x^{tk})^t\over 1-x^k}
\Bigl(
\sum_{\l\in\setP} x^{t|\l|} y^{\#\{h\in\setH(\l), h=1\}} 
\Bigr)^t \cr
&=  \prod_{k\geq 1}{(1-x^{tk})^t\over 1-x^k}
 \Bigl(
\prod_{m\geq 1}{1+(y-1)x^{tm}\over 1-x^{tm}}
\Bigr)^t
 \cr
&=  \prod_{k\geq 1}{(1+(y-1)x^{tk})^t\over 1-x^k}.
\qed \cr
}
$$
\medskip

\def\sec{5}
\section{\sec. Other Specializations} 
Some specializations are given in the introduction and in [Ha08a]. 
In this section we collect other specializations of Theorem 1.3.
When the specialization is easy to derive, a simple comment is written 
between brackets.

\proclaim Corollary \sec.1 [$z=0$]. 
We have 
$$
\sum_{\l\in\setP} x^{|\l|} y^{\#\setH_t(\l)} 
= 
\prod_{k\geq 1}
{
(1-x^{tk})^t
\over
(1-(yx^t)^k)^{t}(1-x^k)  
}.\leqno{(\sec.1)}
$$

\proclaim Corollary \sec.2 [$z=0, y=-1$]. 
We have 
$$
\sum_{\l\in\setP} x^{|\l|} (-1)^{\#\setH_t(\l)} 
= 
\prod_{k\geq 1}{
(1-x^{4tk})^t (1-x^{tk})^{2t}
\over
(1-x^{2tk})^{3t} (1-x^{k})
}.
\leqno{(\sec.2)}
$$

{\it Proof}.
First, we prove (\sec.2) when $t=1$. By Corollary~\sec.1  
$$
\leqalignno{
\sum_\l (-x)^{|\l|} 
&= \prod_k {1\over 1-(-x)^k} \cr
&= 
\prod_{k} {1\over 1-x^{2k}} 
\prod_{k\ {\rm odd}} {1-x^k\over 1-x^{2k}} \cr
&= \prod_{k} {1-x^{4k}\over (1-x^{2k})(1-x^{4k})} 
 \prod_{k\ {\rm odd}} {1-x^k\over 1-x^{2k}} \cr
&= \prod_k{1-x^{4k}\over (1-x^{2k})^2}\prod_{k\ {\rm odd}} {(1-x^k)} \cr
&= \prod_k{(1-x^{4k})(1-x^k)\over (1-x^{2k})^3}.\cr
}
$$
To obtain identity (\sec.2) for arbitrary $t$ we substitute 
$x:=yx^t$ into the above identity.\qed
\medskip

\proclaim Corollary \sec.3 [$y=1$]. 
We have 
$$
\sum_{\l\in\setP} x^{|\l|} \prod_{h\in \setH_t(\l)} 
\bigl(1-{tz\over h^2}\bigr) 
= 
\prod_{k\geq 1}
{
(1-x^{tk})^z
\over
1-x^k
}.\leqno{(\sec.3)}
$$

Corollary \sec.3 can be seen as a discrete interpolation between formulas
(1.10) and (1.11). For example, we have
$$
\leqalignno{
\sum_{\l} x^{|\l|} \prod_{h\in \setH_1(\l)} 
\bigl(1-{36\over h^2}\bigr) 
&= \prod_{k\geq 1} { (1-x^{k})^{36} \over 1-x^k };\cr
\sum_{\l} x^{|\l|} \prod_{h\in \setH_2(\l)} 
\bigl(1-{36\over h^2}\bigr) 
&= \prod_{k\geq 1} { (1-x^{2k})^{18} \over 1-x^k };\cr
\sum_{\l} x^{|\l|} \prod_{h\in \setH_3(\l)} 
\bigl(1-{36\over h^2}\bigr) 
&= \prod_{k\geq 1} { (1-x^{3k})^{12} \over 1-x^k };\cr
\sum_{\l} x^{|\l|} \prod_{h\in \setH_6(\l)} 
\bigl(1-{36\over h^2}\bigr) 
&= \prod_{k\geq 1} { (1-x^{6k})^6 \over 1-x^k },\cr
}
$$
where each sum is over all 6-cores $\l$.

\proclaim Corollary \sec.4 [\hbox{$z=-b/y, y\rightarrow 0$}]. 
We have 
$$
\sum_{\l\in\setP} x^{|\l|} \prod_{h\in \setH_t(\l)} 
{tb\over h^2}
= 
e^{bx^t}
\prod_{k\geq 1}
{
(1-x^{tk})^t
\over
1-x^k
}.\leqno{(\sec.4)}
$$

{\it Proof}. 
Using identity (3.9) 
the right-hand side of (1.4) can be written:
$$
\prod_{k\geq 1} { (1-x^{tk})^t \over (1-(yx^t)^k)^{t}(1-x^k)  }
\exp{\Bigl(-z\sum_{m\geq 1} {(yx^t)^m\over m(1-(yx^t)^m)}\Bigr)}.
\leqno{(\sec.5)}
$$
Since
$$
\leqalignno{
\exp{\Bigl({b\over y}\sum_{m\geq 1} {(yx^t)^m\over m(1-(yx^t)^m)}\Bigr)}
&=\exp{\Bigl({b\over y}\bigl({yx^t\over 1-yx^t} + O(y^2)\bigr)\Bigr)}\cr
&=e^{bx^t}+O(y), \cr
}
$$
we obtain (\sec.4) when $z=-b/y$ and $y\rightarrow 0$ in Theorem 1.3 under 
the form (\sec.5). \qed
\medskip

\proclaim Corollary \sec.5 
[{\rm Compare the coefficients of $b^nx^{tn}$ in (\sec.4)}]. 
We have 
$$
\sum_{\l\vdash tn, \#\setH_t(\l)=n}\quad \prod_{h\in \setH_t(\l)} 
{1\over h^2}
= { 1 \over t^n n!  }.\leqno{(\sec.6)}
$$

Formula (\sec.6) is a classical result (see, e.g., [St99, p.469]).
\goodbreak

\proclaim Corollary \sec.6
[{\rm Compare the coefficients of $b^nx^{tn+m}$ in (\sec.4)} ]. 
We have 
$$
\sum_{\l\vdash tn+m, \#\setH_t(\l)=n}\quad \prod_{h\in \setH_t(\l)} 
{1\over h^2}
= { c_t(m) \over t^n n!  },\leqno{(\sec.7)}
$$
where $c_t(m)$ is the number of $t$-cores of $m$.

\goodbreak

\proclaim Corollary \sec.7 [{\rm  Compare 
the coefficients of $(-z)^{n-1}x^{nt}y^n$}]. We have
$$
\sum_{\l\vdash nt, \#\setH_t(\l)=n } \quad
\prod_{h\in\setH_t(\l)}{1\over h^2}\sum_{h\in\setH_t(\l)} h^2 
= {3n-3+2t\over 2 (n-1)!}.
\leqno{(\sec.8)}
$$

{\it Proof}.
Let $R$ be the right-hand side of (1.4). As $R$ is equal to (\sec.5), we have
\goodbreak
$$
\leqalignno{
&[(-z)^{n-1}x^{nt}y^n] R\cr
=&
[x^{nt}y^n] 
{1\over (n-1)!}\prod_{k\geq 1} { (1-x^{tk})^t \over (1-(yx^t)^k)^{t}(1-x^k)  }
\Bigl(\sum_{m\geq 1} {(yx^t)^m\over m(1-(yx^t)^m)}\Bigr)^{n-1}\cr
=&
[x^{t}y] 
{1\over (n-1)!} { 1 \over (1-yx^t)^{t}  }
\Bigl( {1\over (1-(yx^t))}+{yx^{t}\over 2(1-(yx^t)^2)}\Bigr)^{n-1}\cr
=&
[x^{t}y] 
{1\over (n-1)!}  (1+tyx^t)
\Bigl( {(1+yx^t)}+{yx^{t}\over 2}\Bigr)^{n-1}\cr
=&
{1\over (n-1)!}  \bigl((n-1){3\over 2}+t\bigr).\qed\cr
}
$$
\medskip

The above corollary is the $t$-core analogue of the marked hook formula [Ha08a].
When $t=1$ Formula (\sec.8) reduces to 
$$
\sum_{\l\vdash n} \quad
\prod_{h\in\setH(\l)}{1\over h^2}\sum_{h\in\setH(\l)} h^2 
= {3n-1\over 2 (n-1)!}.
$$
We recover the marked hook formula (Theorem 1.5) 
thanks to the famous hook formula
due to Frame, Robinson and Thrall [FRT54]
$$
f_\l={n!\over \prod_{v\in\l} h_v(\l)}, 
$$
where $f_\l$ is the number of standard Young tableaux of shape $\l$
(see [St99, p.376; Kn98, p.59; Kr99; Ze84; GNW79; GV85; NPS97; RW83]).

\proclaim Corollary \sec.8 [$y=1$; {\rm  compare 
the coefficients of $z$}]. We have
$$
\sum_{\l\in\setP } x^{|\l|} 
\sum_{h\in\setH_t(\l)}{1\over h^2}
= 
{1\over t}\prod_{m\geq 1}{1\over 1-x^m}
\sum_{k\geq 1}{ x^{tk}\over k (1-x^{tk}) } .
\leqno{(\sec.9)}
$$

{\it Proof}.
Let $y=1$. Using (3.9) we have
$$
\sum_{\l\in\setP} x^{|\l|} \prod_{h\in \setH_t(\l)} 
\bigl(1-{tz\over h^2}\bigr) 
= 
\prod_{k\geq 1} { 1 \over (1-x^k)  }
\exp{\Bigl(-z\sum_{m\geq 1} {(x^t)^m\over m(1-(x^t)^m)}\Bigr)}.
$$
Comparing the coefficients of $z$ in the above identity yields (\sec.9).\qed
\medskip

\proclaim Corollary \sec.9.  We have
$$
\sum_{\l\in\setP } x^{|\l|} 
\sum_{h\in\setH(\l), h {\rm\ odd}}{1\over h^2}
= 
\prod_{m\geq 1}{1\over 1-x^m}
\sum_{k\geq 1}{ x^{2k}+2x^k\over 2k (1-x^{2k}) } .
\leqno{(\sec.10)}
$$

{\it Proof}.
Let $t=1$ and $t=2$  in Corollary \sec.8. We obtain respectively
$$
\sum_{\l\in\setP } x^{|\l|} 
\sum_{h\in\setH(\l)}{1\over h^2}
= 
\prod_{m\geq 1}{1\over 1-x^m}
\sum_{k\geq 1}{ x^{k}\over k (1-x^{k}) }. \leqno{(\sec.11)}
$$
$$
\sum_{\l\in\setP } x^{|\l|} 
\sum_{h\in\setH_2(\l)}{1\over h^2}
= 
{1\over 2}\prod_{m\geq 1}{1\over 1-x^m}
\sum_{k\geq 1}{ x^{2k}\over k (1-x^{2k}) }. \leqno{(\sec.12)}
$$
Taking the difference between identities (\sec.11) and (\sec.12) 
yields (\sec.10).\qed 
\smallskip

{\it Remark}. Identity (\sec.11) has a direct proof, which makes use of  
an elegant result on 
multi-sets of hook lengths and multi-sets of partition parts 
obtained by
Stanley, Elder, Bessenrodt, Bacher and Manivel et al.
[Be98, BM02, Ho86, St08, KS82, We1, We2].
See [Ha08a] for more details and applications. 
\def\sec{6}
\section{\sec. Improvement of a result due to Kostant} 
Let 
$$
\prod_{n\geq 1} (1-x^n)^s  =  \sum_{k\geq 0} f_k(s) x^k. \leqno{(\sec.1})
$$
Kostant proved the following result [Ko04, Th. 4.28]. 

\proclaim Theorem \sec.1 [Kostant].
Let $k$ and $m$ be two positive integers such that $m\geq \max(k,4)$. 
Then $f_k(m^2-1)\not=0$. 

The condition $m>1$ in the original statement of Kostant's Theorem 
should be replaced by
$m\geq 4$, as, for example, $f_3(8)=0$ (see Theorem~\sec.2).
Our Theorem 1.6 extends Kostant's result in two directions: first,
we claim that $(-1)^k f_k(s)>0$ instead of $f_k(s)\not=0$; second,
$s$ is any real number instead of an integer of the form $m^2-1$. 
\medskip

\smallskip
{\it Proof of Theorem 1.6}.
By identity (1.3)  we may write 
$$
(-1)^kf_k(s)=\sum_{\l\vdash k} W(\l), 
\leqno{(\sec.2)}
$$
where
$$
W(\l)= \prod_{v\in\l} \Bigl({s+1\over h_v^2}-1\Bigr)
= \prod_{v\in\l} \Bigl({s+1-h_v^2\over h_v^2}\Bigr). 
\leqno{(\sec.3)}
$$
For each $\l\vdash k$ and $v\in\l$ we have $h_v(\l)\leq k$, so that
$W(\l)\geq 0$.  This means that there is {\it no cancellation} in 
the sum  (\sec.2).
If $s>k^2-1$, then $W(\l)>0$. If  $s=k^2-1\geq 15$, we have $k\geq 4$.
In that case there is at least one partition $\l$, whose hook lengths
are strictly less than $k$. Hence $W(\l)>0$. \qed 

\smallskip

Here is another result of Kostant [Ko04, Th.4.27] on which we will make some
comments.
\proclaim Theorem \sec.2 [Kostant].
We have
$$
\leqalignno{
f_4(s)&=1/4!\ s(s-1)(s-3)(s-14); \cr
-f_3(s)&=1/3!\ s(s-1)(s-8);\cr
f_2(s)&=1/2!\ s(s-3).\cr
}
$$

Even though we do not see how to factorize each $f_k(s)$, the occurrences 
of some factors in the above formulas have some relevance in terms of 
hook lengths.
Every partition contains one hook length $h_v=1$, so that $f_k(s)$ has
the factor $s+1-h_v^2=s$ (see (\sec.3)). Every partition of $3$ contains
a hook length $h_v=3$, so that $f_3(s)$ has the factor $s-8$. Every
partition of $2$ or $4$ has a hook length $h_v=2$, so that 
$s-3$ is a factor of $f_2(s)$ and $f_4(s)$. 
Note that Lehmer's conjecture claims that 24 is never a root of $f_k(s)$ 
for any positive integer $k$ (see [Se70]).  

\def\sec{7}
\section{\sec. Reversion of the Euler Product} 
Let $y(x)$ be a formal power series satisfying the following relation
$$
\leqalignno{
x&=y(1-y)(1-y^2)(1-y^3)\cdots &{(\sec.1)}\cr
&=y-y^2-y^3+y^6+y^8-y^{13}-y^{16}+\cdots\cr
}
$$
The first coeficients of the reversion series in (\sec.1)
are the following
$$
y(x)=
x+{x}^{2}+3\,{x}^{3}+10\,{x}^{4}+38\,{x}^{5}+153\,{x}^{6}+646\,{x}^{7}
+\cdots \leqno{(\sec.2)} 
$$
They are referred to as the first values of the sequence A109085
in The On-Line Encyclopedia of Integer Sequences [Slo].

\proclaim Theorem \sec.1.
We have the following explicit formula for the reversion of (\sec.1)
in terms of hook lengths:
$$
y(x)=\sum_{n\geq 1}{x^n\over n} \sum_{\l\vdash n-1}\ \prod_{v\in\l}
\bigl(1+{n-1\over h_v^2} \bigr). \leqno{(\sec.3)}
$$

{\it Proof}. Rewrite (\sec.1) as $y=x\phi(y)$ where 
$\phi(y)=\prod_{m\geq 1} (1-y^m)^{-1}$.
By the Lagrange inversion formula and identity (1.3)  we have
$$
\leqalignno{
[x^n]\ y &= {1\over n} [x^{n-1}]\ \phi(x)^n \cr
&= {1\over n} [x^{n-1}]\ {\prod_{m\geq 1} (1-y^m)^{-n}} \cr
&= {1\over n} [x^{n-1}]\ {\sum_{\l\in\setP} \prod_{v\in\l} 
\bigl( 1+{n-1\over h_v^2}\bigr)x
} \cr
&= {1\over n}
\sum_{\l\vdash n-1}\ \prod_{v\in\l}
\bigl(1+{n-1\over h_v^2} \bigr).\qed \cr
}
$$

{\it Proof of Theorem 1.7}.
The first part of Theorem 1.7 is easy to verify by using identity (1.3).
As the coefficients of $y(x)$ defined by (\sec.1) 
are all positive integers, the above theorem implies that the expression
$$
{1\over n+1}\sum_{\l\vdash n} \prod_{v\in\l} \bigl(1+{n\over h_v^2}\bigr)
$$
is a positive integer. \qed


\bigskip
\medskip
{\bf Acknowledgements}.
The author wishes to thank Dominique Foata for 
helpful discussions during the preparation of this paper.
He also thanks Mihai Cipu, Kathy Ji, Alain Lascoux and
Richard Stanley for their comments,
and Andrei Okounkov for pointing to him that the main identity of [Ha08a] is
in [NO06].

\medskip

\vfill\eject
\bigskip \bigskip


\centerline{References}

{\eightpoint

\bigskip 
\bigskip

\divers AF02|Adin, Ron M.; Frumkin, Avital|
Rim Hook Tableaux and Kostant's $\eta$-Function Coefficients,
{\it arXiv: math.CO/0201003}|

\livre An76|Andrews, George E.|The Theory of
Partitions|Addison-Wesley, Reading, {\oldstyle 1976}
({\sl Encyclopedia of Math. and Its Appl.,} vol.~{\bf 2})|

\article Be98|Bessenrodt, Christine|On hooks of Young diagrams|Ann. of 
Comb.|2|1998|103--110|

\divers BM02|Bacher, Roland; Manivel, Laurent|Hooks and Powers of Parts in 
Partitions, {\sl S\'em.  Lothar. Combin.}, 
vol.~{\bf 47}, article B47d, {\oldstyle 2001}, 11 pages|

\divers BG06|Berkovich, Alexander; Garvan, Frank G.|%
The BG-rank of a partition and its applications, 
{\it arXiv: math/0602362}|

\divers CFP05|Cellini, Paola; Frajria, Pierluigi M.; Papi, Paolo|%
The $\hat W$-orbit of $\rho$, Kostant's formula for powers of the Euler 
product and affine Weyl groups as permutations of $\setZ$, 
{\it arXiv: math.RT/0507610}|

\divers CO08|Carlsson, Erik; Okounkov, Andrei|%
Exts and Vertex Operators,
{\it arXiv:0801. 2565v1 [math.AG]}|

\article Dy72|Dyson, Freeman J.|Missed opportunities|%
Bull. Amer. Math. Soc.|78|1972|635--652|

\divers Eu83|Euler, Leonhard|The expansion of the infinite product 
$(1-x)(1-xx)(1-x^3)(1-x^4)(1-x^5)(1-x^6)$ etc. into a single series, 
{\sl English translation from the Latin by Jordan Bell} 
on {\it arXiv:math.HO/0411454}|

\article FK99|Farkas, Hershel M.; Kra, Irwin|On the Quintuple Product
Identity|Proc. Amer. Math. Soc.|27|1999|771--778|

\divers FH99|Foata, Dominique; Han, Guo-Niu|%
The triple, quintuple and septuple product identities revisited. 
{\sl Sem. Lothar. Combin.} Art. B42o, 12 pp|

\article FRT54|Frame, J. Sutherland; Robinson, Gilbert de Beauregard;        
Thrall, Robert M.|The hook graphs of the symmetric groups|Canadian 
J. Math.|6|1954|316--324|

\article GV85|Gessel, Ira; Viennot, Gerard|%
Binomial determinants, paths, and hook length formulae|Adv. in  
Math.|58|1985|300--321|

\article GKS90|Garvan, Frank; Kim, Dongsu; Stanton, Dennis|Cranks and 
$t$-cores|Invent. Math.|101|1990|1--17|

\article GNW79|Greene, Curtis; Nijenhuis, Albert;
Wilf, Herbert S.|A probabilistic proof of a formula for the number of     
Young tableaux of a given shape|Adv. in Math.|31|1979|104--109|

\divers Ha08a|Han, Guo-Niu|An explicit expansion formula for the powers of 
the Euler Product in terms of partition hook lengths, 
{\sl arXiv:0804.1849v2, Math.CO}, 35 pages, {\oldstyle 2008}|

\divers Ha08b|Han, Guo-Niu|Discovering hook length formulas by 
expansion technique, {\sl in preparation}, 42 pages, {\oldstyle 2008}|

\article Ho86|Hoare, A. Howard M.|An Involution of Blocks in the Partitions 
of $n$|Amer. Math. Monthly|93|1986|475--476|

\article JS89|Joichi, James T.; Stanton, Dennis|An
involution for Jacobi's identity|Discrete Math.|73|1989|261--271|

\livre JK81|James, Gordon; Kerber, Adalbert|The
representation theory of the symmetric group|Encyclopedia
of Mathematics and its Applications, {\bf 16}. Addison-Wesley Publishing, 
Reading, MA, {\oldstyle 1981}|

\article Ka74|Kac, Victor G.|Infinite-dimensional Lie algebras and 
Dedekind's $\eta$-function|Functional Anal. Appl.|8|1974|68--70|

\livre Kn98|Knuth, Donald E.|The Art of Computer Programming,  {\bf vol.~3}, 
Sorting and Searching, 2nd ed.|Addison Wesley Longman,  {\oldstyle 1998}|

\article Ko76|Kostant, Bertram|On Macdonald's $\eta$-function formula, the
Laplacian and generalized exponents|Adv. in Math.|20|1976|179--212|

\article Ko04|Kostant, Bertram|Powers of the Euler product and commutative 
subalgebras of a complex simple Lie algebra|Invent. Math.|158|2004|181--226|

\article Kr99|Krattenthaler, Christian|\ Another involution
principle-free bijective proof of  
Stanley's hook-content formula|J.      
Combin. Theory Ser. A|88|1999|66--92| 

\article KS82|Kirdar, M. S.; Skyrme, Tony H. R.|On an Identity Related to 
Partitions and Repetitions of Parts|Canad. J. Math.|34|1982|194-195|

\livre La01|Lascoux, Alain|Symmetric Functions and Combinatorial Operators on 
Polynomials|CBMS Regional Conference Series in Mathematics, Number 99, 
{\oldstyle 2001}|

\article Ma72|Macdonald, Ian G.|\ Affine root systems and 
Dedekind's $\eta $-function|Invent. Math.|15|1972|91--143|

\livre Ma95|Macdonald, Ian G.|Symmetric Functions and Hall Polynomials|
Second Edition, Clarendon Press, Oxford, {\oldstyle 1995}|

\article Mi85|Milne, Stephen C.|An elementary proof of
the Macdonald identities for $A\sp {(1)}\sb l$|Adv. in Math.|
57|1985|34--70|

\article Mo75|Moody, Robert V.|\ Macdonald identities and Euclidean 
Lie algebras|Proc. Amer. Math. Soc.|48|1975|43--52|

\article NPS97|Novelli, Jean-Christophe; Pak, Igor;
Stoyanovskii, Alexander V.|A direct bijective proof of the hook-length
formula|Discrete Math. Theor. Comput. Sci.|1|1997|53--67|

\divers NO06|Nekrasov, Nikita A.; Okounkov, Andrei|Seiberg-Witten theory and 
random partitions. The unity of  mathematics, 525--596, {\sl Progr. Math.}, 
{\bf 244}, Birkhaeuser Boston. {\oldstyle 2006}. (See also 
{\sl arXiv:hep-th/0306238v2}, 90 pages, {\oldstyle 2003})|

\article RS06|Rosengren, Hjalmar; Schlosser, Michael|%
Elliptic determinant evaluations and the Macdonald identities 
for affine root systems|Compositio Math.|142|2006|937-961|

\article RW83|Remmel, Jeffrey B.; Whitney, Roger|A
bijective proof of the hook formula for the number of column strict tableaux
with bounded entries|European J. Combin.|4|1983|45--63|

\livre Se70|Serre, Jean-Pierre|Cours d'arithm\'etique|Collection SUP: 
``Le Math\'emati\-cien", 2 Presses Universitaires de France, Paris 
{\oldstyle 1970}|

\divers Slo|Sloane, Neil; {\it al.}|
The On-Line Encyclopedia of Integer Sequences, {
\tt http:// www.research.att.com/\char126njas/sequences/}|

\divers St08|Stanley, Richard P.|{\sl Errata and Addenda to Enumerative 
Combinatorics Volume 1, Second Printing}, version of 25 April 2008.
{\tt http://www-math.mit.edu/\char126rstan/ ec/newerr.ps}|

\livre St99|Stanley, Richard P.|Enumerative Combinatorics, vol. 2|
Cambridge university press, {\oldstyle 1999}|

\divers Ve|Verma, Daya-Nand|Review of the paper ``Affine root systems and 
Dedekind's $\eta $-function" written by Macdonald, I. G., MR0357528(50\#9996),
MathSciNet, 7 pages|

\divers We1|Weisstein, Eric W.|Elder's Theorem, from MathWorld -- 
A Wolfram Web Resource|

\divers We2|Weisstein, Eric W.|Stanley's Theorem, from MathWorld -- 
A Wolfram Web Resource|

\article Wi69|Winquist, Lasse|%
An elementary proof of $p(11m+6)\equiv 0\,({\rm mod} 11)$|%
J. Combinatorial Theory|6|1969|56--59|

\article Ze84|Zeilberger, Doron|A short hook-lengths bijection
inspired by the Greene-Nijen\-huis-Wilf proof|Discrete 
Math.|51|1984|101--108|

\bigskip

\irmaaddress
}
\vfill\eject

\end